\documentclass[journal]{IEEEtran}
\newcommand{\subparagraph}{}
\usepackage{titlesec}
\titlespacing{\subsection} {0pt}{6pt}{3pt}

\usepackage{cite}

\ifCLASSINFOpdf
  \usepackage[pdftex]{graphicx}
\else
\fi
\usepackage{amssymb}
\usepackage{amsmath}
\usepackage{algorithmic}

\usepackage{enumitem}

\usepackage{multirow}
\usepackage{array}

\usepackage{color,soul}
\definecolor{Blue}{rgb}{0.3,0.3,0.9}
\DeclareMathSizes{6}{6}{6}{6}
\newcommand{\squeezeup}{\vspace{-2.5mm}}
\begin{document}
%
\title{Real-Time Local Volt/VAR Control Under External Disturbances with High PV Penetration}
%
%
%

\author{A.~Singhal,~\IEEEmembership{Member,~IEEE,}
        V.~Ajjarapu,~\IEEEmembership{Fellow,~IEEE,}
         J.C.~Fuller,~\IEEEmembership{Senior~Member,~IEEE,}
        and~J.~Hansen,~\IEEEmembership{Member,~IEEE}%
\thanks{This work was supported in part by U.S. Department of Energy’s Sunshot Initiative Program DE-0006341.}        
\thanks{A. Singhal and V. Ajjarapu are with Department of Electric and Computer Engineering, Iowa State University, Ames, IA, 50010 USA (e-mail: ankit@iastate.edu, vajjarap@iastate.edu).}
\thanks{J. C. Fuller, J. Hensen are with Pacific Northwest National Laboratory, Richland, WA, 99453 USA (e-mail: jason.fuller@pnnl.gov, jacob.hansen@pnnl.gov).}}

%
%

 \markboth{Accepted in IEEE Transactions of Smart Grid - DOI: https://doi.org/10.1109/TSG.2018.2840965}%
 {}
%



\maketitle


\begin{abstract}
Volt/var control (VVC) of smart PV inverter is becoming one of the most popular solutions to address the voltage challenges associated with high PV penetration. This work focuses on the local droop VVC recommended by the grid integration standards IEEE1547, rule21 and addresses their major challenges i.e. appropriate parameters selection under changing conditions, and the control being vulnerable to instability (or voltage oscillations) and significant steady state error (SSE). This is achieved by proposing a two-layer local real-time adaptive VVC that has two major features i.e. a) it is able to ensure both low SSE and control stability simultaneously without compromising either; and b) it dynamically adapts its parameters to ensure good performance in a wide range of external disturbances such as sudden cloud cover, cloud intermittency, and substation voltage changes. A theoretical analysis and convergence proof of the proposed control is also discussed. The proposed control is implementation friendly as it fits well within the integration standard framework and depends only on the local bus information. The performance is compared with the existing droop VVC methods in several scenarios on a large unbalanced 3-phase feeder with detailed secondary side modeling.
\end{abstract}

\vspace{-2.5mm}
\begin{IEEEkeywords}
solar photovoltaic system, smart grid, volt/var control, smart inverter, real-time control, distributed control.
\end{IEEEkeywords}
\vspace{-2mm}
%
\IEEEpeerreviewmaketitle

\section{Introduction}
%
%
%
%


\IEEEPARstart{S}{olar} photovoltaic (PV) penetration is continuously rising, and is expected to be tripled in the next 5 years in the USA \cite{seia_2016}.
High PV penetration is being fueled by the favorable policies and significant cost reductions, nonetheless, it brings its own set of technical challenges such as voltage rise and rapid voltage fluctuations due to cloud transients which could lead to the reduced power quality \cite{coster_integration_2011,mather_high-penetration_2016}. 
In traditional volt/var control (VVC), voltage regulating devices such as capacitors and load tap changers are supposed to maintain the feeder voltage but they are not fast enough to handle transient nature of solar generation i.e. cloud cover \cite{yeh_adaptive_2012,robbins_two-stage_2013,mcgranaghan_advanced_2008}. Therefore, PV inverter has emerged as an effective VVC solution to handle rapid variations in the modern distribution system by providing faster and continuous control capability in contrast to slower and discrete response of traditional devices \cite{turitsyn_options_2011,bollen_voltage_2005,carvalho_distributed_2008}.

The PV inverter VVC methods primarily fall into two broad categories: 1) optimal power flow (OPF) based centralized and distributed control approaches and 2) local control approaches. Most of the literature deals with the OPF based methods which are solved either in a centralized manner \cite{xu_multi-timescale_2017,dallanese_optimal_2014,yeh_adaptive_2012,farivar_optimal_2012} or using distributed algorithms \cite{zheng_fully_2016,zhang_optimal_2015,chen_robust_2017,robbins_two-stage_2013,turitsyn_distributed_2010}. There are several other distributed control methods proposed for PV inverter VVC which can be referred from the latest comprehensive survey papers \cite{antoniadou-plytaria_distributed_2017} and \cite{molzahn_survey_2017}.  However, the extensive communication requirements among the PV devices challenge the real-time implementation of these methods. Additionally, communication delays and the large time requirement to solve most OPFs limit their ability to respond to faster disturbances at seconds time scale such as cloud intermittency \cite{turitsyn_options_2011,carvalho_distributed_2008,zhu_fast_2016}. Though distributed algorithms are relatively faster, most of these methods assume constant substation voltage and rely on full feeder topology information for control parameter selection which is usually not fully known to the utilities or not always reliable. These issues make OPF based VVC methods difficult to implement and also vulnerable to fast external disturbances such as cloud transients, changes in substation voltage and topology changes. Therefore, we focus on the local VVC approaches in this work which are usually faster, simple to implement, and can respond to the sudden external disturbances in the distribution systems.

Among local approaches, droop VVC is the most popular local control framework among utilities and in the existing literature. It was first proposed by \cite{seal_standard_2010} which now has been adopted by the IEEE1547 integration standard \cite{noauthor_ieee_2018} and also being widely used by Rule 21 in California \cite{noauthor_rule_nodate}. Local control is simple to implement based on the local bus information, however, ensuring the system-wide control stability and performance is a challenge in the local control design. It has been identified an improper selection of control parameters can lead to control instability and voltage oscillation issues \cite{andren_stability_2015,farivar_equilibrium_2013,jahangiri_distributed_2013}. Most literature on droop control \cite{zhang_three-stage_2017,malekpour_dynamic_2017,karthikeyan_coordinated_2017} or other similar local control methods \cite{shah_online_2016,safavizadeh_voltage_2017} lack in analytical characterization and do not discuss the parameter selection and the control stability/convergence issues. Some work such as delayed droop control in \cite{jahangiri_distributed_2013}  discuss stability issue and scaled var control in \cite{zhu_fast_2016} provide rigorous performance analysis. But none of them adapt themselves in changing operating conditions and external disturbances to ensure control convergence. Another major challenge with droop control is its inherent inability to achieve a low steady state error (SSE) while ensuring convergence in all conditions. In other words, a certain slope selection which ensures the control stability, may also lead to high SSE as indicated by \cite{farivar_equilibrium_2013}; and as shown later in this paper, both the control stability and low SSE are crucial for the distribution systems operations. 

In this work, our focus is to analyze and design a droop-based local VVC which addresses following two challenges associated with the conventional droop VVC: a) To make the control parameters selection self-adaptive to changing operating conditions and external disturbances; and b) To achieve both low SSE and control stability simultaneously without compromising either. To achieve these objectives, we propose a fully local and real-time adaptive VVC within the IEEE1547 standard framework and provides its convergence and performance (SSE) analysis. We also compare our proposed VVC with the conventional droop VVC and it's improved version 'delayed' droop VVC \cite{jahangiri_distributed_2013} which improves the stability performance under normal condition but is vulnerable to control instability and high SSE under external disturbances due to lack of proper parameter selection, as detailed soon.
Our work extends the previous works and provides unique contributions in following way: 1) The proposed control achieves both low SSE and control stability simultaneously by decoupling the two objectives; 2) The control parameters are made self-adaptive to  commonly occurring external disturbances such as cloud intermittency, cloud cover, changing load profile, and substation voltage changes; 3) A theoretical analysis of convergence of the proposed adaptive VVC is discussed and a sufficient condition for convergence is derived; 4) It is compatible with the existing onboard droop controls specified in the recent standard (IEEE1547); 5) The real-time adaptive nature and tight voltage control feature of the proposed control opens interesting opportunities for operators to utilize PV inverters not only to mitigate over-voltage but for other volt/var related applications such as CVR, loss minimization, providing var support to transmission side etc.; and 6) A detailed modeling of secondary side of an unbalanced distribution system is used to verify the control approach with house-level loads and heterogeneous inverter population.

The layout of the remainder of the paper is as follows. In Section II, stability conditions and SSE expressions of the conventional droop VVC are derived and discussed to establish the base for the adaptive VVC development. Based on the analysis, the adaptive control strategy is developed in Section III. Section IV provides the convergence analysis of the proposed VVC with an illustration. Simulation results on the test system are discussed in Section V. Finally, concluding remarks are presented in Section VI.

\begin{figure}[t]
	\centering
	\includegraphics[trim=0in 0in 0in 0in,width=2.2in]{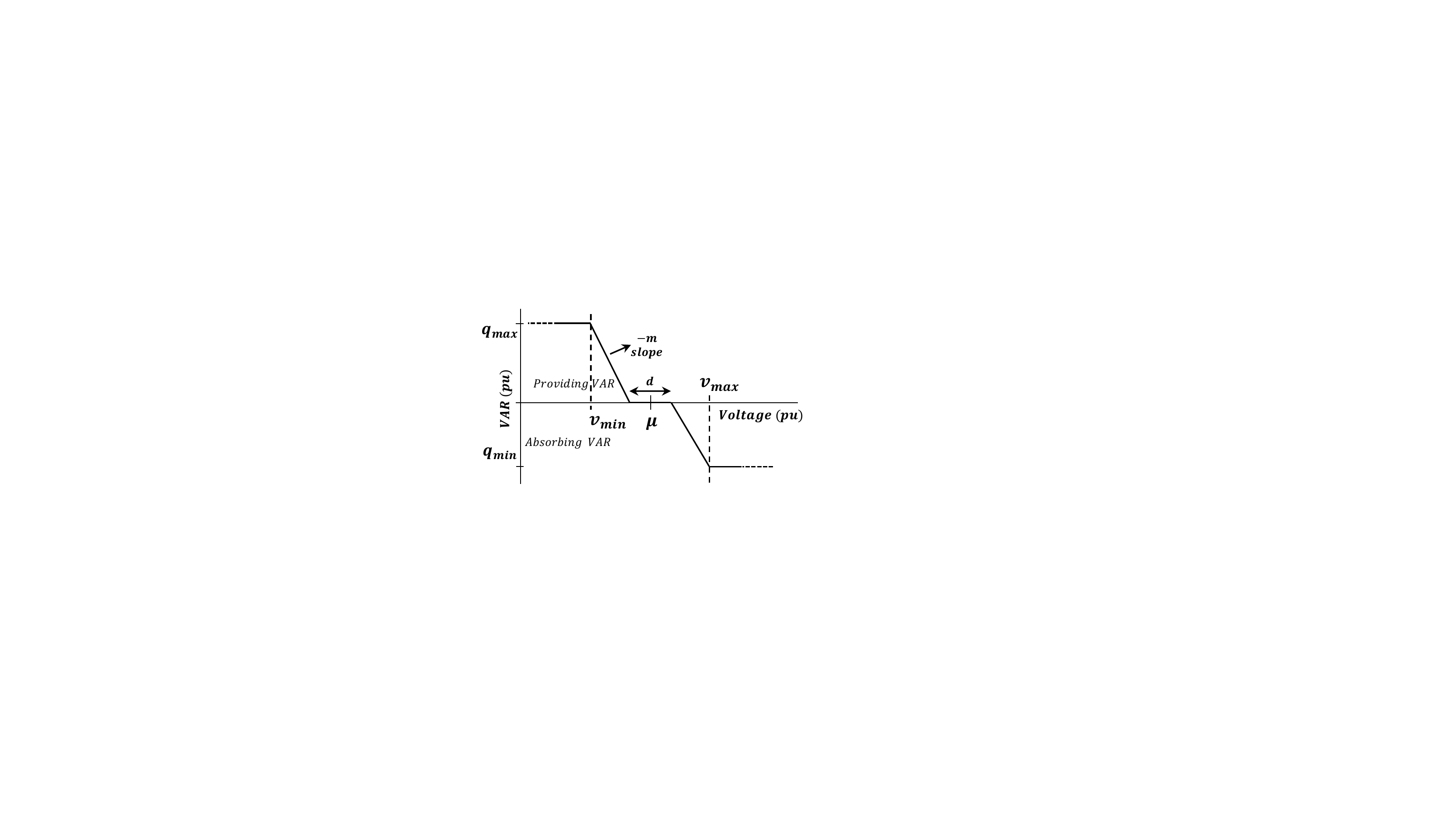}
    \squeezeup
    \caption {Conventional droop VVC framework recommended by IEEE 1547 }
    \vspace{-2mm}
    \label{fig:std_droop}
\end{figure}

\section{Background And Problem Setup }
Consider a general $N+1$ bus distribution system with one substation bus and $N$ load buses with PV inverters. The power flow equations for the system can be written as
\begin{equation}
\label{eq:power_flow}
\begin{split}
P^{inv}-P_d= g_p(V,\delta); \quad Q^{inv}-Q_d= g_q(V,\delta)\\
\end{split}
\end{equation}
Where 
$Q\!=\![Q_2 Q_3 \dots Q_{N+1}]^T$ and $P^{inv}\!=\![P_2 P_3 \dots P_{N+1}]^T$ are inverter reactive and real power injection vectors respectively at each bus. $P_d$ and $Q_d$ are similar vectors of real and reactive power loads at each bus. $g_p$ and $g_q$ are well-known power flow equations with voltage magnitude ($V$) and angles ($\delta$) as variables at all load buses \cite{grainger_power_1994}. 
The standard droop function $f_i (.)$ at $i^{th}$ bus is shown in \figurename\ref{fig:std_droop}. It is a piecewise linear function with a deadband $d$ and slope $m$. Assuming the operating point is in non-saturation region, the inverter var dispatch at time $t$ can be written as a function of previous voltage and other control parameters as
\begin{equation}
\label{eq:droop}
{{Q}_{i,{t+1}}}=f_i(V_{i,t})=-m_i({{V}_{i,t}}-\mu_i \pm d/2)
\end{equation}
where, $Q_{i,t}$ and $V_{i,t}$ are the inverter var injection and the voltage magnitude respectively. The sign $\pm d/2$ represents that $-d/2$ is used for $v_i>\mu_i$ and $+d/2$ is used for $v_i<\mu_i$. Subscripts  $i$ and $t$ denote $i^{th}$ bus and time instant $t$. $\mu_i$ is the reference voltage and $m_i$ is the slope of the curve.  We consider the same slope for both the regions in the droop control for a given inverter as shown in \figurename \ref{fig:std_droop}. $m_i$ can be maintained at desired value by changing control parameters as
\begin{equation*}
\label{eq:slope}
m_i=q_{i,max}/(\mu_i-d/2-v_{i,min})= q_{i,min}/(\mu_i+d/2-v_{i,max})
\end{equation*}
where $v_{i,min}$, $v_{i,max}$, $q_{i,min}$, $q_{i,max}$ are the four control set-points. It should be noted that, in the existing droop methods, these parameters are either constant or un-controlled. Whereas, in this work, these parameters are dispatched based on the proposed adaptive control strategy. As detailed soon, dynamic control over these parameters leads to more reliable control performance compared to the previous works.
\subsection{Stability Analysis}
As described in \cite{farivar_equilibrium_2013,jahangiri_distributed_2013},  the local droop VVC can be modeled as feedback dynamical system $\phi$ with $N$ states $[Q_{2,t}\;Q_{3,t}\; \dots\;Q_{N+1,t}]^T$ at discrete time $t$.
\begin{equation}
\label{eq:phi}
Q_{t+1}=\phi(Q_t)=f(h(Q_t))
\end{equation}
Where the vector $f(.)\!=\![f_2\;f_3\;\dots\;f_{N+1}]$ contains local VVC functions which map the current voltage vector $V_t$ to new inverter var injections vector $Q_{t+1}$ i.e. $Q_{i,t+1}\!=\!f_i (V_{i,t})$. The new var vector $Q_{t+1}$, in turn, leads to the new voltage vector $V_{t+1}$ according to power flow equations (1). The function $h$ is an implicit function vector derived from (1) i.e. $h_i (Q_{i,t} )=V_{i,t}$. It is shown in \cite{jahangiri_distributed_2013} that the system $\phi$ is locally stable in the vicinity of an equilibrium point ($\bar{Q}$) if all eigenvalues of the matrix $\partial \phi/\partial Q$ have magnitude less than 1.
\begin{equation}
\small
\label{eq:phi_slope}
\left[\frac{\partial\phi}{\partial Q}\right]_{Q=\bar{Q}}=
\left[\frac{\partial f}{\partial V}\right]
\left[\frac{\partial V}{\partial Q}\right]
\end{equation}
In the case of droop control, $\partial f/\partial V$ is a diagonal matrix with slope at each inverter as diagonal entries. 
\begin{equation}
\small
\label{eq:M}
\left[\frac{\partial f}{\partial V}\right]=-M=-diag(m_i)=-
\begin{bmatrix}
m_2 & \cdots & 0\\ 
\vdots & \ddots & \vdots\\
0 & \cdots & m_{N+1}\\ 
\end{bmatrix}
\end{equation}
Let’s define {\small $A\!=\!\partial V/\partial Q$} and  {\small $a_{ij}\!=\!\partial V_i/\partial Q_j$} which is a voltage sensitivity matrix with respect to var injection and can be extracted from the power flow Jacobian matrix from (1) as shown in \cite{jahangiri_distributed_2013}.   
\begin{equation}
\small
\label{eq:A}
\left[\Delta V\right]=\left[A\right]\left[\Delta Q\right]=
\begin{bmatrix}
a_{22} & \cdots & a_{2,N+1}\\ 
\vdots & \ddots & \vdots\\
a_{N+1,2} & \cdots & a_{N+1,N+1}\\ 
\end{bmatrix} \left[\Delta Q\right]
\end{equation} 
In other words, the sufficient condition for the control stability can be written as
\vspace{-2mm}
\begin{equation}
\small
\label{eq:rho}
\rho(MA)<1
\vspace{-3mm}
\end{equation}
Where $\rho$ is the spectral radius of a matrix which is defined as the largest absolute value of its eigenvalues. Condition (\ref{eq:rho}) provides useful information for evaluating the stability of specific inverter slope settings. However, in order to obtain information for selecting the inverter slopes, we will derive another conservative sufficient condition for stability using spectral radius upper bound theorem \cite{horn_matrix_2012}.

\textit{Theorem 1:} Let $\|.\|$ be any matrix norm on $\mathbb{R}^{n \times n}$ and let $\rho$ be the spectral radius of a matrix, then for all $X \in \mathbb{R}^{n \times n}$:
\begin{equation}
\label{eq:rho_X}
\rho(X)\leq ||X||
\end{equation}
\textit{Proposition:} If sum of each row of $MA$ is less than 1, i.e.
\begin{equation}
\label{eq:slope_cond}
m_i. \sum_{j=1}^{N}|a_{ij}|< 1 \quad \forall i ,
\end{equation}
Then the droop control will be stable i.e. $\rho(MA)<1$

\textit{Proof:} Using Theorem 1, if we apply {\small $\|.\|_\infty$} on {\small $MA$}, then,
{\small $\rho(MA)\leq \|MA\|_\infty=\max_{1\leq i\leq N}\sum_{j=1}^{N}|m_i.a_{ij}|$}. If condition (\ref{eq:slope_cond}) holds true i.e. {\small $m_i. \sum_{j=1}^{N}|a_{ij}|< 1 \quad \forall i$}, then the maximum of sum of rows will also be less than one. Thus, the upper bound on spectral radius will always be less than one i.e. {\small $\rho(MA)<1$}.  

\textit{Remark 1:} The condition (\ref{eq:slope_cond}) provides useful information for slope selection for each inverter to ensure control stability, i.e. {\small $m_i\!<\!m_i^c$}, where $m_i^c$ is critical slope given by 
$m_i^c\!=\!(\sum_j|a_{ij}|)^{-1}$
It should be noted that, usually, the entries of the sensitivity matrix $A$ do not remain constant. Changes in operating conditions (cloud cover, load changes) as well as changes in feeder topology lead to change in values of $a_{ij}$ and $m_i^c$; thus, they can potentially cause instability, if $m_i$ are not updated dynamically. Intuitively, entries of $A$ can also be seen as proportional to the reactance of the feeder lines \cite{farivar_equilibrium_2013} i.e. longer lines are more likely to have higher magnitude of $a_{ij}$ and lower value of critical slope. Therefore, PV inverters on rural network with longer lines, especially towards the feeder end, will be more sensitive to instability and their slope selection should be more conservative. Therefore, non-adaptive and homogeneous slope selection for all inverters make system prone to control instability. 
It is worth mentioning that an attempt to lower the effective slope by adding a delay block after droop in the delayed droop \cite{jahangiri_distributed_2013} improves the stability compared to the conventional droop i.e. $Q_{t+1}=f_i(V_{i,t})+\tau. Q_t$, where $\tau$ is a delay coefficient. However, because of its non-adaptive nature and un-controlled parameters, it may lead to issues under external disturbances and topology changes which will be illustrated through a comparison later in this section.

\subsection{Steady State Error (SSE) Concerns}
One of the major drawbacks of the droop control is the significant deviation from the set-point in steady state. To derive the analytical expression for SSE, let’s assume the system is at equilibrium point ($\overline{Q},\overline{V}$) at $t=0$. For simplicity, let's also assume $d=0$ in this analysis (it can be extended for non-zero $d$ values also, but this will make the analysis unnecessarily complicated, and will distract the reader from the main purpose of the section which is to illustrate the SSE concerns of the conventional droop controls). Control equation (2) can be written in vector form at $t=0$, as
\begin{equation}
\small
\label{eq:eqb}
[\overline{Q}]=-[M][\overline{V}-\mu]
\end{equation}
Now, consider an external disturbance perturb the equilibrium by causing sudden change in the voltage, $\Delta V^d$, at $t=0$ which changes the voltage at $t=0$ i.e. {\small $[V]_{t=0}=\overline V+ \Delta V^d$}. This drives control to dispatch the new var at $t=1$ i.e. $Q_{t=1}=-[M][V_{t=0}-\mu]$. Now, the following can be written,
\begin{equation}
\small
\label{eq:eqb}
[Q_{t=1}-\overline{Q}]=[\Delta Q]_{t=1}=-[M]\Delta V^d
\end{equation}
Using the similar procedure, following can be written for $t>0$
\begin{equation}
\small
\label{eq:diff1}
[\Delta Q]_{t+1}=-[M][\Delta V]_t
\end{equation}
Where {\small $[\Delta Q]_{t+1}=[Q_{t+1}-Q_t]$} and {\small $[\Delta V]_t=[V_t-V_{t-1}]$}.Using (6) and (12), we can write
\begin{equation}
\small
\label{eq:diff2}
[\Delta V]_{t+1}=-[A][M][\Delta V]_t
\end{equation}
\begin{equation}
\small
\label{eq:diff3}
[\Delta V]_{t+1}=[-A.M]^t[\Delta V]_{t=1}
\end{equation}
Using (6) and (11), replace {\small $[\Delta V]_{t=1}=A[\Delta Q]_{t=1}=-[AM]\Delta V^d$} in (14),
\begin{equation}
\small
\label{eq:diff3}
[V]_{t+1}=[V]_t+[-A.M]^{t+1}\Delta V^d
\end{equation}
By writing the (\ref{eq:diff3}) recursively and replacing {\small$[V]_{t=1}=(\overline V+\Delta V^d)-[AM]\Delta V^d$}, 
\begin{equation}
\small 
\label{eq:diff4}
[V]_{t+1}=\overline V+\sum_{i=0}^{t+1}[-A.M]^{i}\Delta V^d
\vspace{-1mm}
\end{equation}
In this case, the geometric progression series of matrices only converges if the condition (\ref{eq:rho}) holds true (the stable case). The new equilibrium voltage can be written as
\begin{equation}
\small
\label{eq:lim}
\lim_{t\to \infty}[V]_{t+1}=\overline V+[I+A.M]^{-1}\Delta V^d
\end{equation}
Finally, SSE vector can be written as
\begin{equation}
\small
\label{eq:sse}
SSE = \lim_{t\to \infty}[V]_{t+1} - \mu
\end{equation}
Equation (\ref{eq:lim}) and (\ref{eq:sse}) show that for a given disturbance, the only way to decrease SSE is to set higher values of slopes $m_i$ which in turn may violate stability condition (\ref{eq:slope_cond}). Usually, SSE is compromised to ensure control stability.

\textit{Remark 2:} Note that, in some cases, it might be possible to maintain voltages within the American National Standard Institute (ANSI) allowable limits with high SSE, close to the boundaries, for a given system condition. But, any external disturbance can instantly push the voltages out of the limits as illustrated later. Moreover, other than complying with ANSI standard, the tight voltage regulation capability (low SSE) makes the system more flexible and provides extra room to the operator to perform other voltage-dependent applications such as CVR, loss minimization etc.; thus fully utilizing the PV inverter’s capability. It can also be shown easily that the delayed droop has the same SSE as the conventional droop.

\begin{figure} [t]
	\centering
	\includegraphics[trim=0in 0in 0in 0in,width=2.6in]{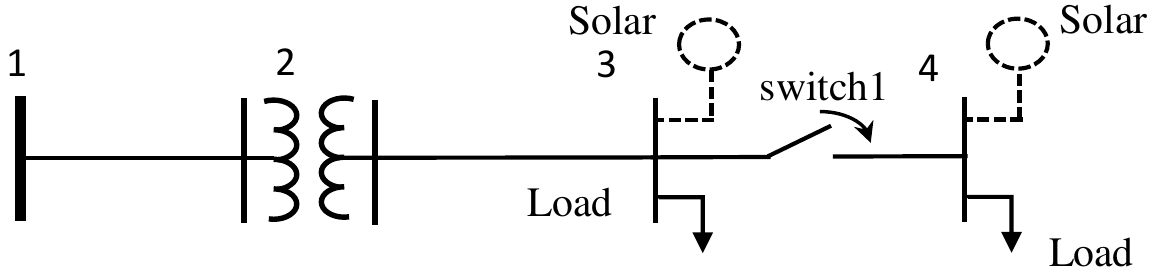}
    \vspace{-1mm}
    \caption {A small 4 bus system to illustrate the impact of external disturbances }
    \label{fig:toy_ckt}
\vspace{2mm}

	\centering
	\includegraphics[trim=0.0in 0in 0in 0in,clip,width=3.3in]{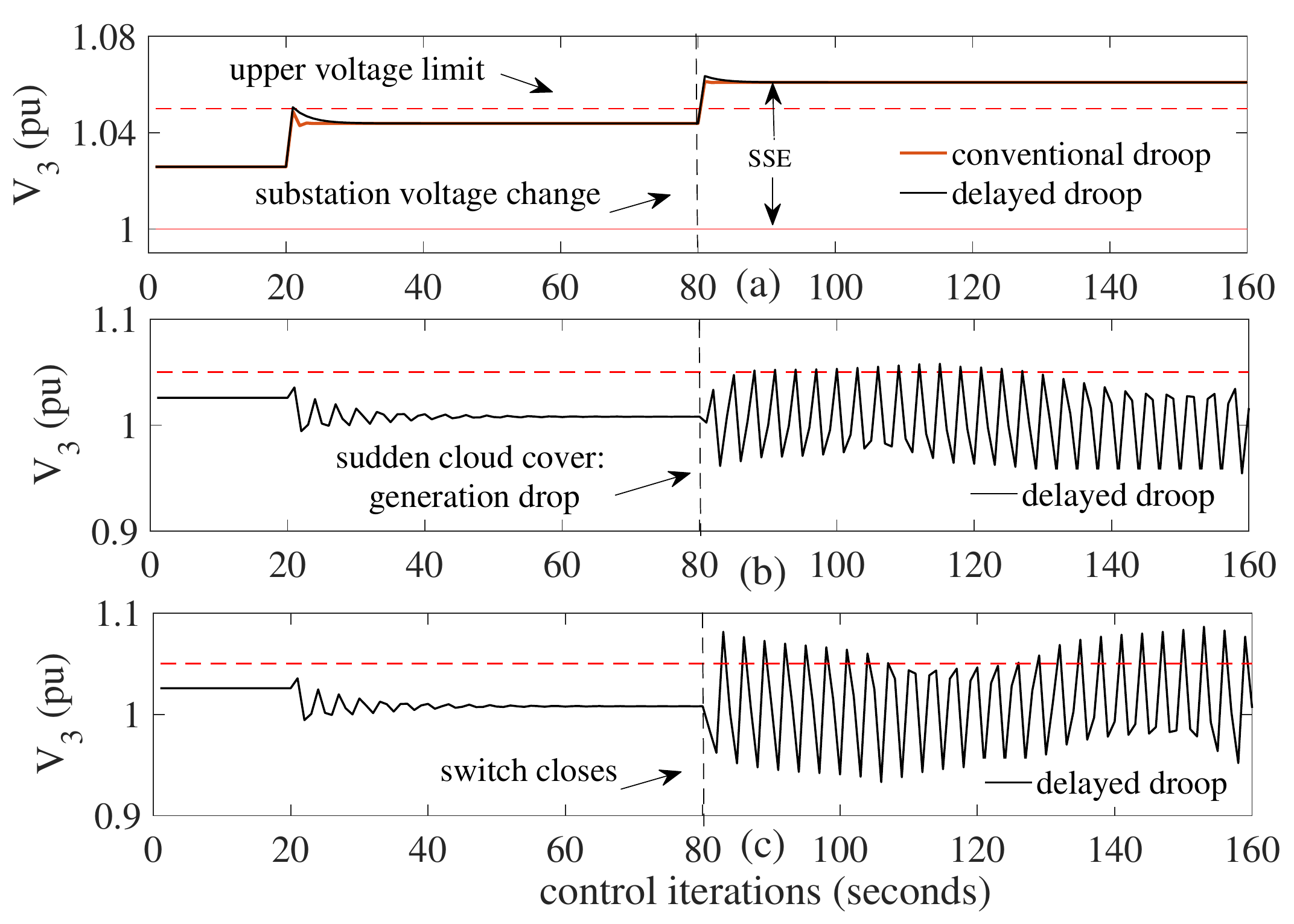}
    \vspace{-2.5mm}
    \caption {Non-adaptive droop VVC performance under impact of a) substation voltage change at conservative slope setting; b) sudden cloud cover at non-conservative setting and; c) topology change at non-conservative setting.}
    \label{fig:toy_results}   
    \vspace{-0.5mm}
\end{figure}

\subsection{Illustration}
To corroborate the above analysis, we will illustrate the impact of external disturbances using a small modified IEEE 4 bus test system shown in \figurename\ref{fig:toy_ckt}. 600 kW load and 900 kW solar generation is added at node 3. A similar node 4 is added via a normally open switch to simulate the change in feeder topology. We will consider two types of initial slope settings to convey the main outcome of the analysis i.e. conservative ($m=1$) and non-conservative ($m=6$). Solar generation is applied at $t=20$ to observe the impact of VVC with $\mu=1$ at node 3 voltage profile. 
\figurename\ref{fig:toy_results}(a) demonstrates how conservative setting causes high SSE (though, within the ANSI limit initially) for the droop controls (both conventional and delayed droop) which leads to over-voltage violation due to a small change in substation voltage from 1.03-1.05 pu at $t=80$. On the other hand, using non-conservative setting to reduce SSE makes the system prone to control instability or voltage flicker as shown in \figurename\ref{fig:toy_results}(b) and (c). Conventional droop is not shown as it is always unstable in these cases. \figurename\ref{fig:toy_results}(b) shows that a sudden drop in solar generation due to cloud cover at $t=80$ increases $q_{max}$ and makes the slope very high which causes voltage oscillations. Further, to simulate the impact of topology change or error in topology information, switch1 is closed at $t=80$. Delayed droop, as discussed before, is stable under normal conditions, however, change in feeder topology leads to voltage oscillations as shown in \figurename\ref{fig:toy_results}(c) at non-conservative settings. This example demonstrates that it is difficult to achieve both low SSE and control stability under external disturbances with the existing droop controllers. Moreover, this problem becomes more crucial in a large realistic system due to thousands of independent inverter devices, higher possibility of inaccuracy in topology information and in parameter selection, and increasing disturbances in this new environment of pro-active distribution system.

Therefore, our intention is to develop a new droop based adaptive VVC strategy 1) to achieve both low SSE and low voltage oscillations (stability) simultaneously; 2) to make control parameters dynamically self-adaptive to external disturbances in real-time.

\vspace{-2mm}
\section{Adaptive Control Strategy}
\begin{figure}[h]
	\centering
	\includegraphics[trim=0in 0in 0in 0in,width=3.42in]{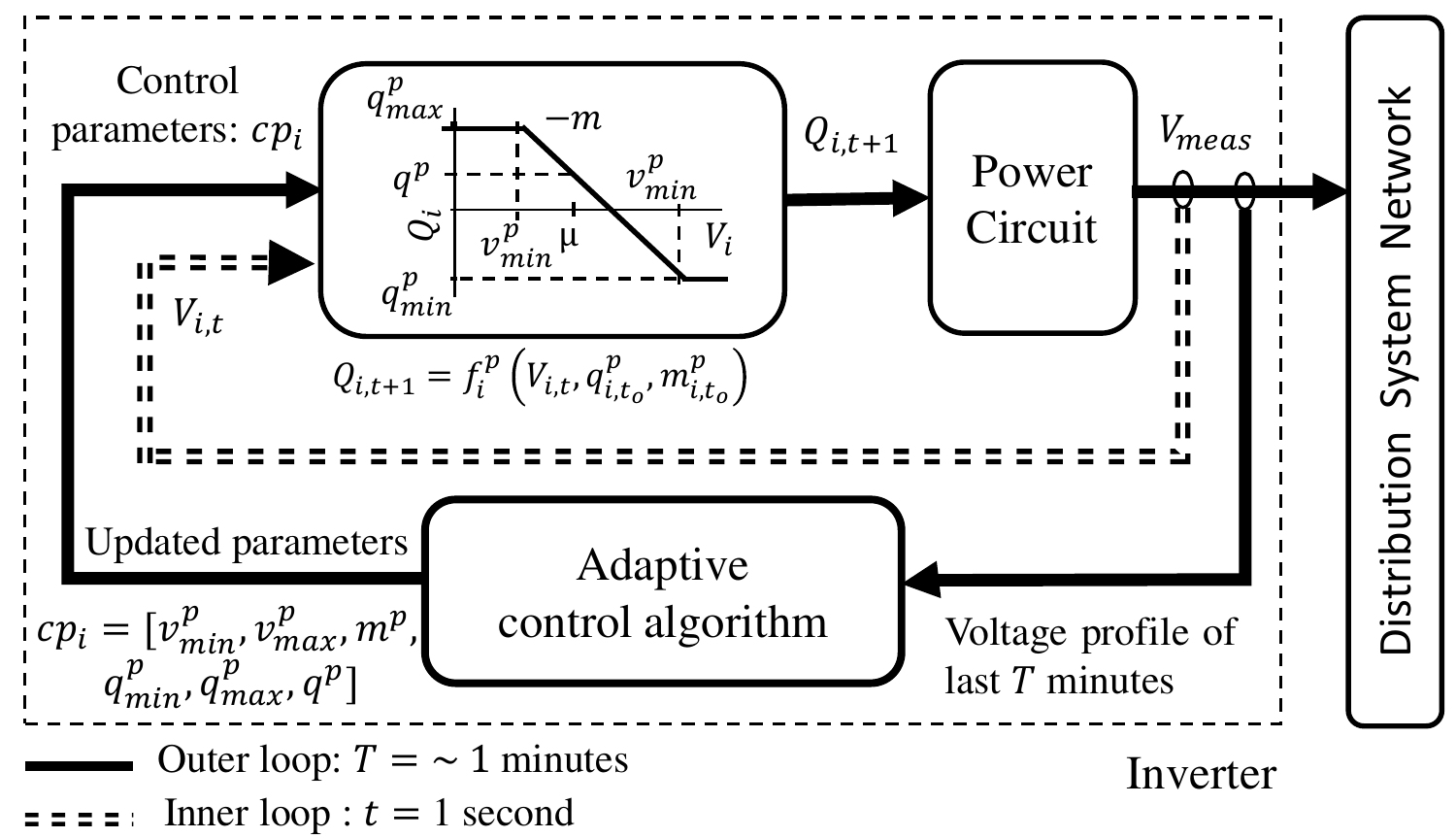}
    \vspace{-1mm}
    \caption {Two-layer framework of the proposed adaptive control approach }
    \vspace{2mm}
    \label{fig:framework}
	\centering
	\includegraphics[trim=0in 0in 0in 0in,width=2.6in]{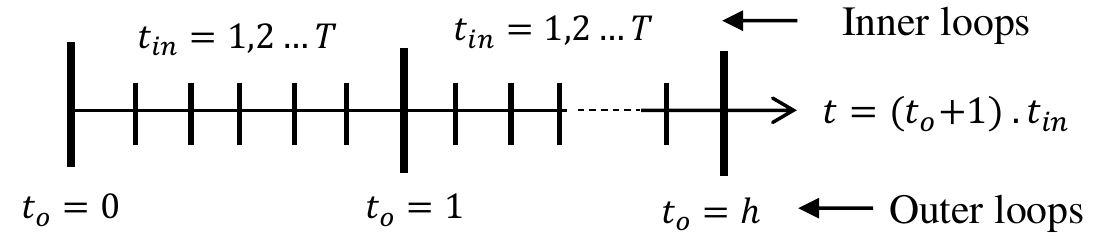}
    \squeezeup
    \caption {Time-line of adaptive inner and outer loop control }
    \label{fig:timeline}
    \vspace{-3mm}
\end{figure}

This section will introduce the proposed adaptive local VVC function $f_i^p (V_{i,t},cp_i)$ which can be written as follows:
\setlength\abovedisplayskip{6pt}
\setlength\belowdisplayskip{6pt}
\begin{equation}
\small
\label{eq:adaptive_droop}
{{Q}_{i,{t+1}}}=f^p_i(V_{i,t},cp_i)=\mathbb{P}[q^p_{i}-m^p_i({{V}_{i,t}}-\mu_i)]
\end{equation}
Where {\small $cp_i=[m_i^p,q_i^p,q_{min,i}^p,q_{max,i}^p,v_{max,i}^p,v_{min,i}^p,]$} are control parameters. The Function {\small $\mathbb{P}$} is a saturation operator with {\small$(q_{min,i}^p,q_{max,i}^p)$} as saturation var limit parameters applied at cut-off parameters {\small $(v_{max,i}^p,v_{min,i}^p)$}. Variable $q_i^p$ is an error adaptive parameter and its main function is to provide SSE correction. Desired adaptive slope $m_i^p$ can be set as,
\begin{equation}
\small
\label{eq:adaptive_slope}
m_i^p=\frac{q_{min,i}^p-q_i^p}{\mu_i-v_{max,i}^p}=\frac{q_{max,i}^p-q_i^p}{\mu_i-v_{min,i}^p}
\end{equation}
There are two unique features of this control. First, the functions of maintaining control stability and low SEE are decoupled. Two different parameters $m^p$ and $q^p$ are used to achieve control stability and low SSE respectively with different approaches so that none of the objectives are compromised. Secondly, all these parameters are dynamically adapted in real-time. Superscript $p$ denotes the adaptive nature of the control parameters. To achieve this, a two-layer control framework is proposed as shown in \figurename \ref{fig:framework}. The inner layer is a fast VVC function $f_i^p (V_{i,t})$ to track the desired set-point $\mu_i$ according to (\ref{eq:adaptive_droop}).  The outer layer dispatches the control parameters $(cp_i)$ based on the proposed adaptive algorithm described later in the section. The outer layer works on a relatively slower time scale $(t_o)$ to allow inner fast control to reach steady state before dispatching new control parameters, thus avoiding hunting and over-corrections. Control time-line is shown in \figurename \ref{fig:timeline}. Control parameters are updated at every period $T$, control horizon of the outer loop control. Each iteration of the inner and outer loop control is denoted by $t_{in}$ and $t_o$ respectively. 
	The adaptive algorithm consists of two strategies where $q_i^p$ and $m_i^p$ are dynamically adapted to take care of SSE and voltage instability/flicker respectively as described below.

\subsection{Error Adaptive Control: Strategy I}
\begin{figure}
	\centering
	\includegraphics[trim=0in 0in 0in 0in,width=2.2in]{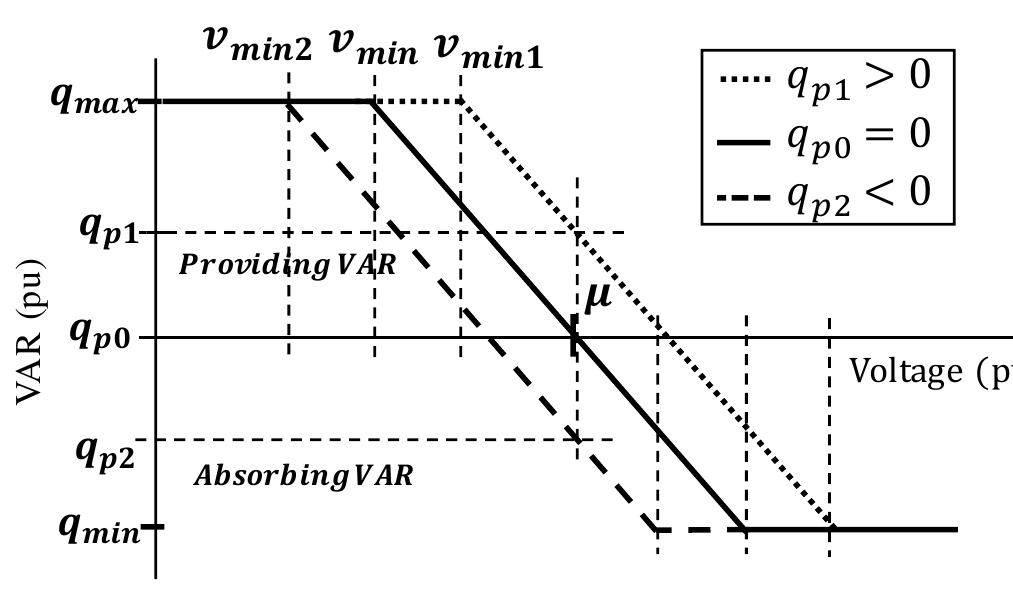}
    \caption {Adaptive strategy I: VVC droop curve with adaptive error adaptive parameter $q_p$ }
    \label{fig:strategy2}
    \squeezeup
\end{figure}

The aim of the strategy I is to minimize SSE by adapting the error adaptive parameter $q_i^p$. We will analyze how the proposed control (\ref{eq:adaptive_droop}) helps to mitigate the SSE and accordingly develop a mechanism to adapt $q_i^p$ locally. Consider the system is at equilibrium point {\small$(\overline{Q},\overline{V})$} at $t=0$ with {\small $SSE=\overline{V}-\mu$}. Now if the parameter $q^p$ is changed at $t=0$ by $\Delta q^p$, the new voltage deviation ($SSE^{adp}$) can readily be obtained by following the procedure provided in the Section II.B  by replacing the conventional droop (\ref{eq:droop}) with the adaptive control (\ref{eq:adaptive_droop}):
\begin{equation}
\small
\label{eq:sse_adp}
SSE^{adp}=\overline{V}-\mu+[I+AM^p]^{-1}A\Delta q^p
\end{equation}
To achieve $SSE^{adp}=0$, $\Delta q^p$ required will be,
\begin{equation}
\small
\label{eq:sse_adp2}
\Delta q_{req}^p=-(A^{-1}+M^p)SSE
\end{equation}	
Equation (\ref{eq:sse_adp2}) provides the analytical expression of the required change in $q^p$ parameter to achieve zero SSE in just one iteration. However, this solution requires the information of $A$ matrix, SSE and slope $(M)$ at all inverter buses which is not available to local bus controllers. Moreover, estimation of $A$ is contingent to error in centralized feeder topology information and may not be reliable. Therefore, we propose a local version of the analytical solution (\ref{eq:sse_adp2}) i.e. {\small $\Delta q_i^p\!=\!-k_i^d.SSE_{avg,i}(t_o)$}, where {\small$SSE_{avg}$} defined for each outer loop as
\setlength\abovedisplayskip{2pt}
\setlength\belowdisplayskip{4pt}
\begin{equation}
\small
\label{eq:sse_avg}
 SSE_{avg,i}(t_0)=\sum_{t_{in}=1}^{T}({V_{t_0t_{in},i}-\mu_i})/{T}
\end{equation}
$SSE_{avg,i}$ denotes the average set-point deviation of voltage at $i^{th}$ inverter bus. A tolerance band for $SSE_{avg,i}$ can be defined as $\mu_i\pm\epsilon_{sse}$, where $\epsilon_{sse}$ is tolerance for the deviation. In this strategy, the adaptive term $q_i^p (t_o)$ in (\ref{eq:adaptive_droop}) is updated at each outer loop interval $t_o$, based on the real-time estimation of $SSE_{avg,i}$ during the last time horizon $T$ as
\setlength\abovedisplayskip{4pt}
\setlength\belowdisplayskip{4pt}
\begin{equation}
\small
\label{eq:qp_update}
q^p_i (t_o )=q^p_i (t_o-1)-k_i^d.SSE_{avg,i}(t_o) 
\end{equation}
Since the $SSE_{avg,i}$ requires only local voltages and only needs to be calculated in the outer loop, there is enough room to calculate this variable without causing any extra delay in the control. It is important to note that {\small $SSE_{avg}$} is used as an algebraic value with sign. The sign of the error decides whether $q_i^p$ needs to be moved positive or negative. If the voltage settles on a higher value than the set point, a negative term is added in $q_i^p$ to facilitate more var absorption to lower the voltage. Similarly, a positive term is added in $q_i^p$ to provide more var when voltage settles lower than the set point. A constant {\small $k_i^d\!>\!0$} is a correction factor which can be decided once from the offline studies.
It's selection  affects the convergence speed of the control which is discussed in detail in Section IV.
\figurename \ref{fig:strategy2} depicts the adaptive control $f_i^p (V_{i,t},cp_i)$ with different $q_i^p$ values. Note that the solid curve with {$q_i^p\!=\!0$} is same as the conventional droop control $f_i (V_{i,t})$ in (2). \figurename \ref{fig:strategy2} brings out an important feature of the proposed control that it can be seen as ``shifted and adaptive" droop VVC which makes it compatible with the integration standards. 

{Nonetheless, it should be noted that the proposed approach may take more than one iterations to achieve near zero SSE, unlike the analytical solution. However, it is compensated by the advantage that it requires only local bus information, thus, making it more feasible. Nevertheless, the update strategy can always be made faster and more accurate using} (\ref{eq:sse_adp2}),{ if information at other nodes is also available in the future. A more detailed theoretical convergence analysis of this adaption strategy is discussed in Section IV.} 

\subsection{Adaptive Slope Control: Strategy II}
The objectives of the strategy II are to ensure stability as well as to keep voltage fluctuations within the IEEE 141 standard\cite{ieee_141_1994} limit by adapting parameter $m^p_i$. Therefore, based on \cite{ieee_141_1994}, we use short-term voltage flicker (VF) as control criteria which is calculated for each inverter bus at the beginning of each outer loop as
\setlength\abovedisplayskip{2pt}
\setlength\belowdisplayskip{2pt}
\begin{equation}
\small
\label{eq:vf}
VF(t_o)=\sum_{t_i=1}^{T}\frac{(V_{t_ot_i}- V_{t_ot_i-1})/ V_{t_0t_i}}{T}\times 100
\end{equation}

\begin{figure}[]
	\centering
	\includegraphics[trim=0in 0in 0in 0in,width=2.2in]{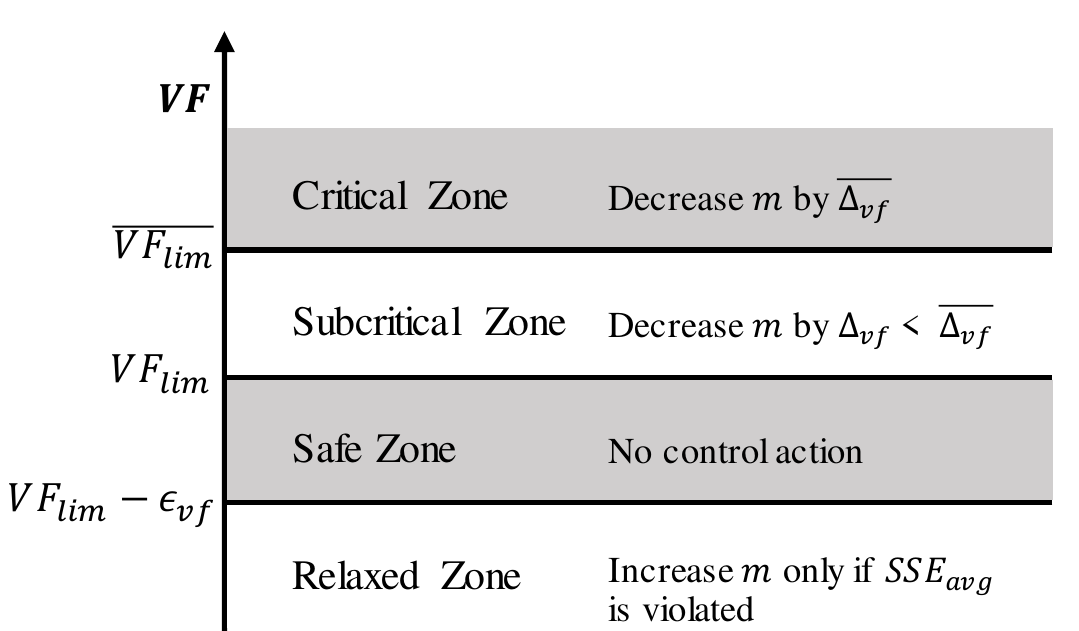}
    \caption {Control action region for adaptive outer loop control strategy II for flicker mitigation }
    \label{fig:flicekr_region}
    \vspace{-2mm}
\end{figure}

As seen in (15), voltage fluctuations are proportional to slope and can be reduced by decreasing $m_i$. For this purpose, the voltage flicker range is divided into four control regions as shown in \figurename \ref{fig:flicekr_region}. The IEEE 141 flicker curve provides the maximum fluctuation limit {\small $(\overline{VF_{lim}})$} beyond which we define as \textit{critical flicker zone}. The same standard also gives a borderline flicker limit {\small $(VF_{lim})$}. The region between {\small$(\overline{VF_{lim}})$} and $VF_{lim}$ is termed as the \textit{subcritical flicker zone}. Further we define a tolerance {\small $(VF_{lim}\!-\!\epsilon_{vf})$} and the tolerance band is termed as the \textit{safe flicker zone}. The region below safe flicker zone is defined as the \textit{relaxed flicker zone}. In critical zone, we update the parameters by a larger amount {\small$(\overline{\Delta_{vf}})$} to avoid control instability and to return to subcritical zone faster. In subcritical zone, the slope is decreased in a smaller step {\small$(\Delta_{vf})$} to avoid over-correction which might impact  {\small$SSE_{avg}$}  negatively. As soon as we reach the safe zone, no control action is taken. This is the desired range of control parameters. Though rarely required, in the relaxed zone, slope is increased to improve SSE only if SSE is out of range. Correction factors {\small$(\Delta_{vf})$} are estimated offline in this work based on the sensitivity analysis and engineering judgment. The amount of slope change required to reduce $VF$ from {\small$(\overline{VF_{lim}})$} to {\small$({VF_{lim}})$} can be taken as an approximate value of {\small$\Delta_{vf}$}. Almost twice of  {\small$\Delta_{vf}$} can be taken as {\small$\overline{\Delta_{vf}}$}. They can also be made responsive to the online control performance, if required.   
\begin{figure}[]
	\centering
	\includegraphics[trim=0in 0in 0in 0in,width=2in]{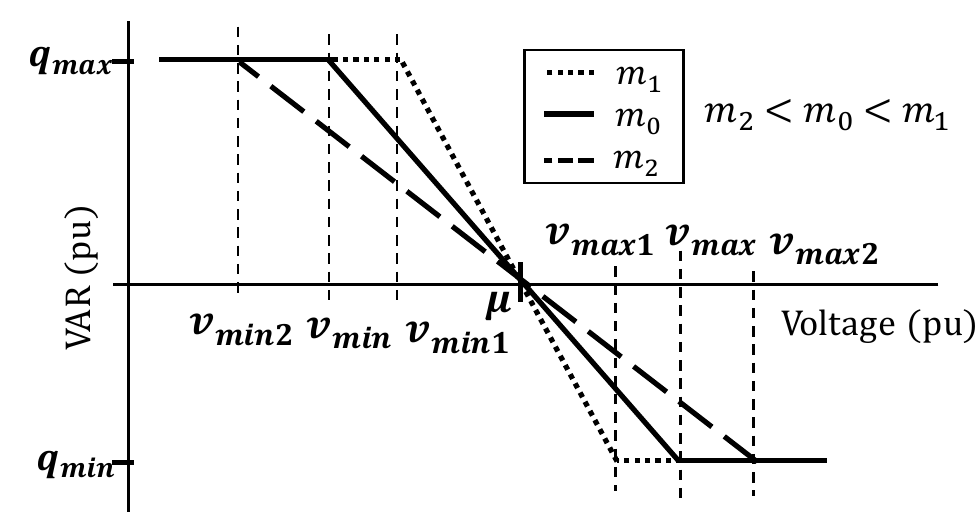}
    \caption {Adaptive strategy II: changing slope of droop curve by changing $v_{min}$ and $v_{max}$ parameters to keep flicker in the limit }
    \label{fig:strategy1}
    \squeezeup
\end{figure}

It is worth noting here that the the main feature of this strategy lies in the decoupling of the two functionalities i.e. SSE and slope. Since SSE is catered by $q_i^p$, slope can always be in the conservative range (safe or relaxed zones) to ensure control stability. {This might cause a momentarily high SSE until the strategy I adapts the SSE again but it prevents the possibility of control instability over large range of operating points}. In this work, we use the earlier derived condition (\ref{eq:slope_cond}) to choose initial slopes. It is estimated using offline studies for the base case, however, to keep safe margin it can be further reduced by a certain factor. \figurename\ref{fig:strategy1}. depicts the control strategy II with adaptive $m^p_i$.

While it is possible to defer real power solar generation, in this work the consumer value is maximized by limiting var output to leftover capacity and not curtailing real power generation. To utilize the inverter capacity entirely, $q_{max}^p$ and $q_{min}^p$ are also updated in every outer loop as 
\begin{equation}
\small
\label{eq:inv_capacity}
q_{max}^p(t_o)\!=\!\sqrt{s^2\!-\!p_{pv}^2(t_o)}; q_{min}^p(t_o)\!=\!-\!\sqrt{s^2\!-\!p_{pv}^2(t_o)}
\end{equation}
where, $s$ is inverter rating and $p_{pv}(t_o)$ is the average solar PV real power generation in the last outer loop time interval. {Note that ideally we would want to use $p_{pv}$ forecast for the next outer loop. Various established methods for short-term solar PV forecasting can be used for this purpose such as statistical or neural network based approaches} \cite{pvforecast_2015}.{ However, since this is not the focus of this paper, for simplicity we use the simplest prediction available which is the previous $p_{pv}$ value.}   

Thus, we get the new parameters {\small $q_i^p (t_o)$} from strategy I, {\small $m_i^p(t_o)$} from strategy II and {\small $ q_{min,i}^p(t_o), q_{max,i}^p(t_o)$} from (\ref{eq:inv_capacity}). 
Finally, {\small $v_{min,i}^p(t_o)$  and $v_{max,i}^p(t_o)$} parameters are calculated using (\ref{eq:adaptive_slope}) 
and dispatched to be used in the inner loop. Overall detailed algorithm of the adaptive control is shown in \figurename \ref{fig:algorithm}. 

\begin{figure}[]
\vspace{-2mm}
\small
\setlist[enumerate,1]{wide=\parindent}
\setlist[enumerate,2]{leftmargin=4.5em}
\setlist[enumerate,3]{leftmargin=4em}
\renewcommand\labelenumi{\theenumi.}
\renewcommand\labelenumii{\theenumi.\arabic{enumii}.}
\renewcommand\labelenumiii{}
\begin{tabular}{p{0.45\textwidth}}
\hline\\[-1.5ex]
\textbf{Algorithm 1}: Adaptive control scheme\\
\hline\\\vspace{-1.5em}
\begin{enumerate}
\item Real-time measurement and control criterion calculation
	\begin{enumerate}
	\item Collect $V_{t=t_o.t_{in}} \forall 			t_{in}=1,2,\dots,n$
	\item Calculate $SSE_{avg}(t_o)$ and $VF(t_o)$
	\end{enumerate}
\item Go to adaptive strategy I: error adaptive   
	\begin{enumerate}
	\item 	If $|SSE_{avg}(t_o)|>\mu+\epsilon_{sse}$
    	\begin{enumerate}
    	\item $q_p(t_o)=q_p(t_o-1)-k_d.SSE_{avg}(t_o)$
	  	\end{enumerate}
	\item Else, $q_p (t_o )=q_p (t_o-1)$
    \end{enumerate}    
\item Go to adaptive strategy II: slope adaptive

\begin{center}
$m^p(t_0)= m^p{(t_0-1)}+\Delta_{m}$
\end{center}
	\begin{enumerate}
	\item If $VF(t_o) > \overline{VF_{lim}}$ \hspace{56pt}	
		$\Delta_m=-\overline{\Delta_{vf}}$
	\item Else if $VF (t_0 )>VF_{lim}$   \hspace{37pt}	
		$\Delta_m=-{\Delta_{vf}}$
    \item Else if $VF (t_0 )>(VF_{lim}-\epsilon_{vf})$  \hspace{6pt}	
		$\Delta_m=0$
	\item Else, check if $|SSE_{avg}|>\mu+\epsilon_{sse}$ \hspace{-2pt}	
		$\Delta_m={\Delta_{vf}}$
	\end{enumerate}
\item Update $q_{max}(t_0)$ and $q_{min}(t_0)$: equation (\ref{eq:inv_capacity})
\item Update final parameters $v_{min}$ and $v_{max}$: equation (\ref{eq:adaptive_slope})
\item  $t_o=t_o+1$, go to step 1
\vspace{-1em}
\end{enumerate}\\
\hline 
\end{tabular}
\caption{Overall algorithm of the proposed adaptive control strategy}
\label{fig:algorithm}
\end{figure}
\section{{Convergence of the Proposed Local Adaptive Control Algorithm}}
The convergence properties and conditions of the proposed adaptive control (\ref{eq:adaptive_droop}) will be investigated in this section. Since it's a two-layer control, we need to study the convergence of both the control loops. If we assume the time horizon $T$ is sufficient for faster control to reach its steady state, the inner loop control equation within the time $T$ can be written as
\vspace{1em}
\begin{equation}
\label{eq:inner_loop}
{{Q}_{i,{t_{in}+1}}}=q^p_{i}-m^p_i({{V}_{i,t_{in}}}-\mu_i)
\end{equation}
Where adaptive parameters $q_i^p$ and $m_i^p$ are constant for the time horizon $T$. Using the stability analysis performed in Section II, a sufficient convergence criteria for (\ref{eq:inner_loop}) can be derived which is same as given by condition (\ref{eq:slope_cond}) and remark 1, i.e. $m_i^p<m_i^c$. Therefore, the inner local control will always converge if the chosen slopes are below critical slope values. Adaptive control strategy II helps to maintain this condition by keeping slope in conservative range as discussed earlier.

\subsection{Outer Loop Control Convergence}
To derive the analytical expression for the outer loop convergence criteria, lets assume the system is currently at outer loop iteration $t_o$ for which inner control has already reached its steady state ($\overline{V_{t_o}}$,$\overline{Q_{t_o}}$). This can be represented as
\vspace{2mm}
\begin{equation}
\small
\label{eq:outer_loop1}
{[\overline{Q}]_{t_o}}=[q^p]_{t_o}-[M^p] [S]_{t_o}
\end{equation}
Where, let's represent the $SSE_{t_0}$ vector with a shorter notation $ S_{t_o}=[\overline{V}_{t_o}-\mu]$. Now the $q^p$ parameter is locally updated for outer loop iteration $(t_o+1)$ based on the update formula given in (\ref{eq:qp_update}) i.e. $q^p_{i,{t_o+1}}=q^p_{i, t_o}-k_i^d.S_{i,(t_o)}$. Now a new updated inner loop steady state ($\overline{V}_{t_o+1}$,$\overline{Q}_{t_o+1}$) is reached for the outer loop iteration $(t_o+1)$ which can be written in vector form as

\begin{equation}
\small
\label{eq:outer_loop2}
{[\overline{Q}]_{t_o+1}}=[q^p]_{t_o}-[K^d][S]_{t_o}-[M^p] [S]_{t_o+1}
\end{equation}
Where $K^d$ is a diagonal matrix with the correction factor $k_i^d$ at each bus as its diagonal entries. Since the objective is to minimize SSE, we can replace all other state variable with $S_{t_o}$ by subtracting (\ref{eq:outer_loop2}) from (\ref{eq:outer_loop1}) and using (\ref{eq:A}) i.e. $\Delta V = A \Delta Q$ as following:
\begin{equation}
\label{eq:sse_conv}
[S]_{t_o}-[S]_{t_o+1}= A\Big[K^d.S_{t_o}-M^p(S_{t_o}-S_{t_o+1})\Big]
\end{equation}
After manipulating (\ref{eq:sse_conv}), the outer loop control can be written as a discrete feedback dynamical system with SSE as the only state variable as following.
\begin{equation}
\label{eq:sse_conv2}
[S]_{t_o+1}= B[S]_{t_o}
\end{equation}
Where $B=I-[I+AM^p]^{-1}AK^d$ is a constant matrix which defines the convergence behavior of the outer loop control. We know that any linear discrete system $x[k+1]=Bx[k]$ is asymptotically stable if and only if all eigenvalues of $B$ have magnitude less than 1 \cite{chen_linear_2013}. Consequently, the convergence condition for the outer loop control can be written as 
\begin{equation}
\label{eq:conv_cond}
\rho(B)<1
\end{equation}
Where $\rho$ is spectral radius of a matrix. Therefore, the convergence of the outer loop is ensured as long as the selection of $K^d$ matrix does not violate the condition (\ref{eq:conv_cond}). In such cases, the outer loop control system (\ref{eq:sse_conv2}) will always converge to zero SSE for any amount of initial SSE i.e. $\lim_{t_o\to\infty} S_{t_o}=0$ for any $S_{t_o=0}$.

\textit{Remark 3:}  {It is interesting to observe here if we choose a non-diagonal $K^d=A^{-1}+M^p$, $B$ turns out to be zero and the SSE can be converged to zero in just one iteration. However, in that case, the $q^p$ update ($q^p_{t_0+1}=q^p_{t_0}-K^d.S_{t_o}$) does not remain local i.e. to update $q^p_i$ at $i^{th}$ node, non-diagonal entries $k^d_{ij}$ need to be multiplied with SSE at all other $j^{th}$ nodes. Therefore, we compromise with the convergence speed to take advantage of the local feature of the control.}

\begin{figure}[h]
\vspace{-0.5em}
	\centering
	\includegraphics[trim=0.0in 0in 0in 0in,width=3.5in]{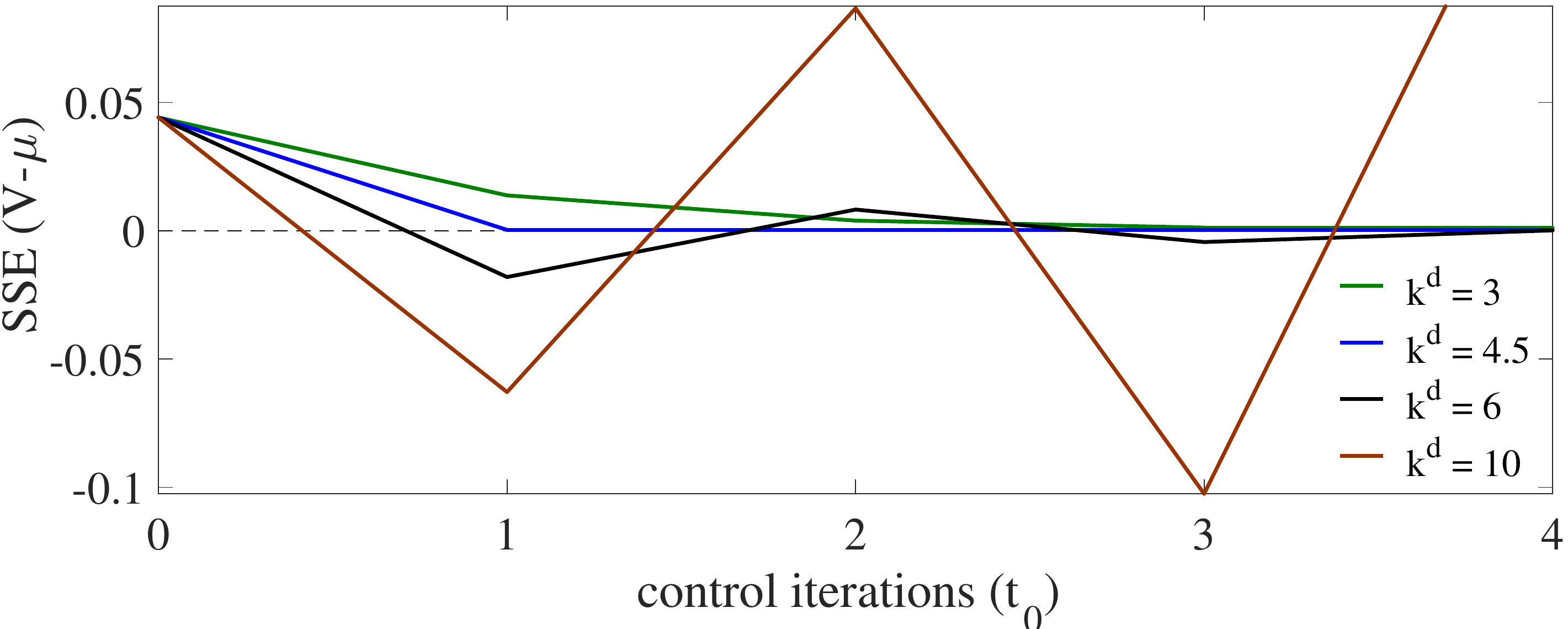}
    \caption {SSE convergence profile of the proposed adaptive outer loop control under different values of $k^d$}
    \label{fig:2bus_kd_sse}   
\end{figure}

For a better understanding of this analysis, let's consider a simple example system described in \figurename \ref{fig:toy_ckt}. Note that this example system can be seen as an equivalent two-bus system as there is only one load bus (bus 3) when switch is open. In that case, $B=b$ and $K=k$ will be a scalars, $b=1-k/(a^{-1}+m)$. The sensitivity matrix $A=a_{33}=0.2857$ var pu/volt pu is calculated offline for the base case. Using 
$m_i^c=(\sum_j{|a_{ij}|})^{-1}$, derived from condition (\ref{eq:slope_cond}), the critical slope for this system turns out to be $m^c=1/a_{33}=3.5$. In order to be in conservative range, $m=1$ is chosen as initial slope. This system converges to zero only if $|b<1|$ i.e. $0<k<2(a^{-1}+m)$. Further, within this stable region, three special cases can be analyzed i.e. $k<(a^{-1}+m), k=(a^{-1}+m)$ and $(a^{-1}+m)<k<2(a^{-1}+m)$. \figurename \ref{fig:2bus_kd_sse} demonstrates the SSE response of the example system under these three stable cases and one unstable case. In the first case of $k^d<4.5$, the system converges to zero without oscillations (overdamped response). In the second case of $k^d=4.5$, system reaches zero SSE in just one iteration and in the third case of $k^d>4.5$, it converges with decaying oscillations (underdamped response). For $k^d>9$, the SSE starts diverging with non-decaying oscillations.

\vspace{-0.5em}
 \section{Case Studies and Discussion}
\subsection{Small System Illustration}
\begin{figure}
	\centering
	\includegraphics[trim=0.0in 0in 0in 0in,width=3.5in]{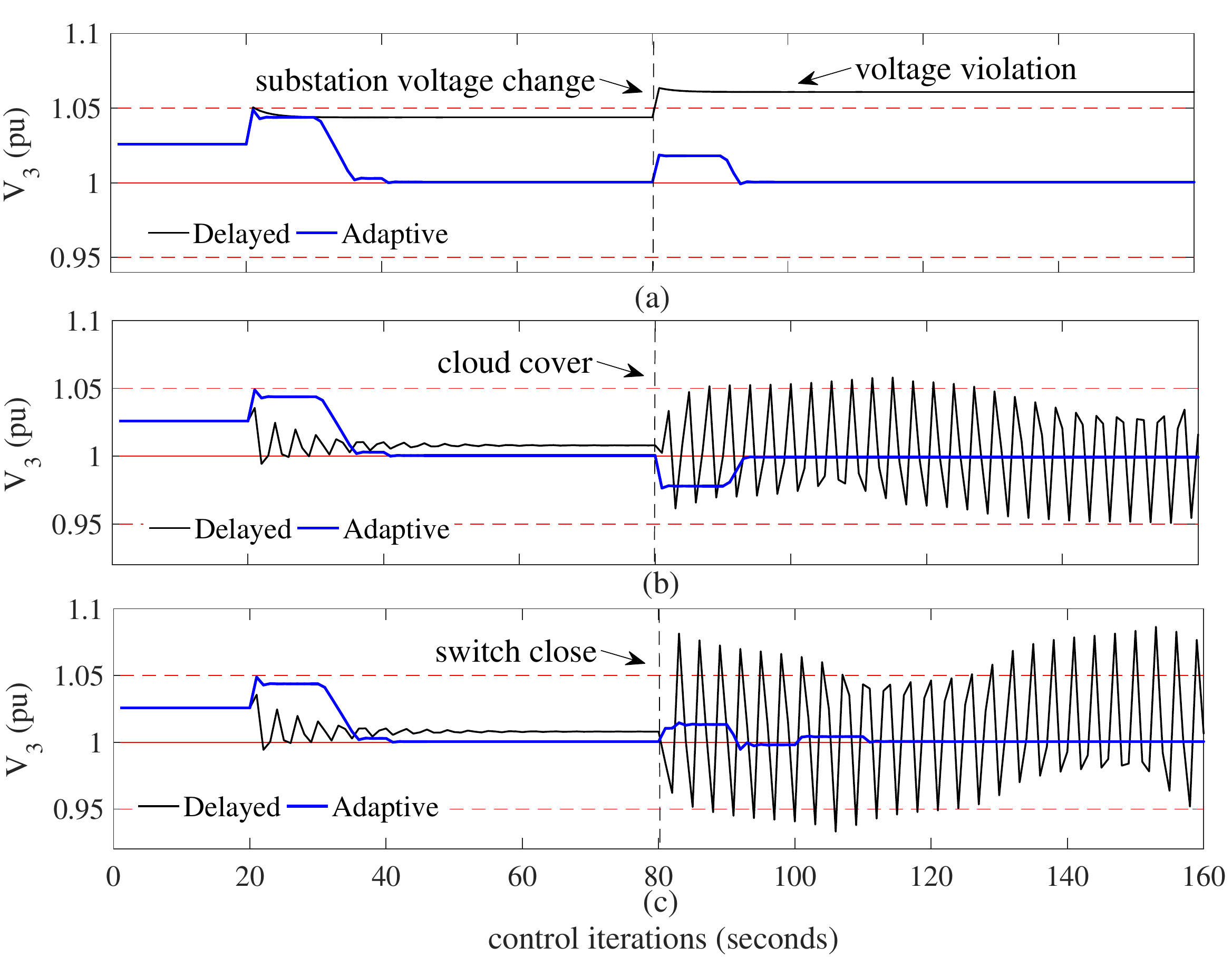}
    \caption { Adaptive VVC performance comparison with delayed VVC under impact of : a) substation voltage change at conservative slope setting; b) sudden cloud cover at non-conservative slope setting and; c) topology change at non-conservative slope setting}
    \label{fig:toy_adaptive}   
\end{figure}

{In this section, the proposed adaptive control performance is discussed and compared with the non-adaptive delayed control for the small example system described in the Section II.C in} \figurename \ref{fig:toy_ckt}. {Since it is a small system, outer loop time horizon of 10 seconds is adequate to demonstrate adaptive nature of the control. Other system setup and parameters selection are same as described earlier. Based on the convergence discussion in the last section, $k^d=4$ and $m=1$ are chosen. }\figurename \ref{fig:toy_adaptive} {compares the performance under three types of disturbances. In all cases, solar generation with VVC is applied at $t=20$. However, it can be seen that the adaptive VVC starts reducing SSE only after 10 seconds as outer adaptive loop works after the time horizon $T$. At $t=80$, when voltage profile gets a surge due to change in substation voltage from 1.03 to 1.05, the adaptive control adapts itself to re-track the set-point within just 1 iterations as visible in }\figurename \ref{fig:toy_adaptive}(a). {Whereas, the delayed VVC leads to voltage violation due to high SSE and non-adaptive nature. Under non-conservative settings, in} \figurename \ref{fig:toy_adaptive} {(b), adaptive control is able to maintain the smooth voltage profile under the impact of sudden cloud cover after $t=80$, unlike the delayed control. Note that for first 10 seconds (from $t=20$ to $t=30$), the adaptive VVC has higher SSE as it uses conservative slope setting and SSE correction starts only after 10 seconds. In case of }\figurename \ref{fig:toy_adaptive}{(c), switch is closed at $t=80$ which reduces the critical slope $m^c$ of the overall system and leads to voltage oscillations in the delayed control. However, the adaptive VVC re-tracks the set-point at both the nodes within few iterations. The new B matrix and convergence condition for $k^d=4$ can be verified as following.}
\begin{equation*}
B=\begin{bmatrix}
0.224 & -0.623\\
-0.646& -0.055\\
\end{bmatrix}, eigenvalues(B)=0.73, 0.56
\end{equation*}

\subsection{Large Test Case Modeling}
\begin{figure}
	\centering
	\includegraphics[trim=0in 0in 0in 0in,width=3.2in]{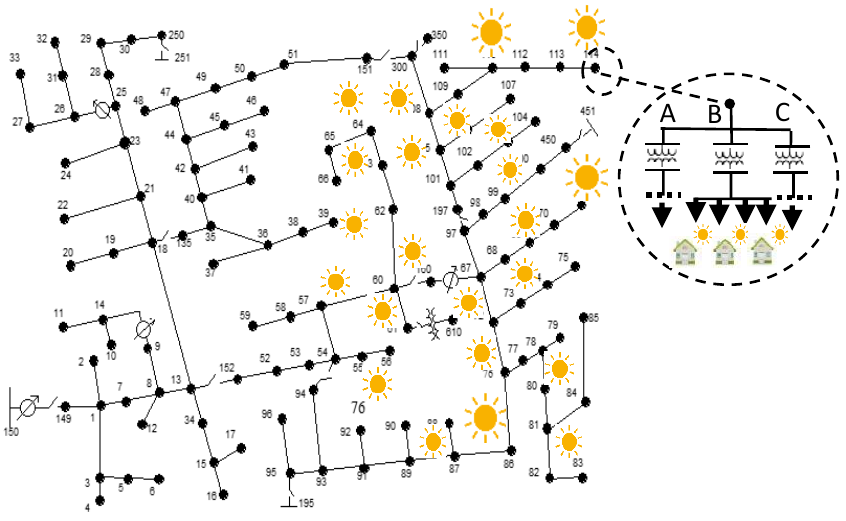}
    \caption {IEEE 123 bus test system with detailed secondary side modeling}
    \label{fig:test_ckt}
    \squeezeup
\end{figure}

\begin{figure}
	\centering
	\includegraphics[trim=0in 0in 0in 0in,width=2.8in]{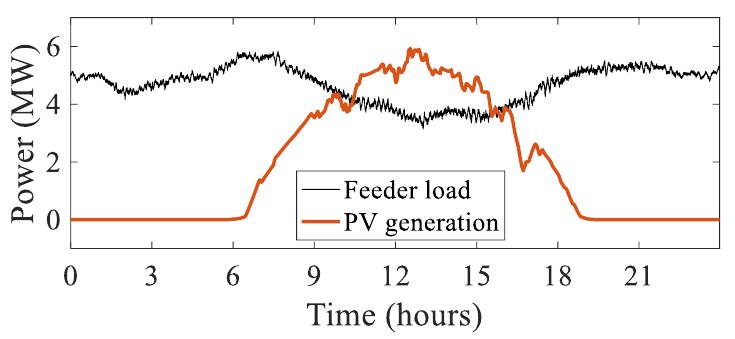}
    \squeezeup
    \caption {Total feeder load and solar PV profile for 24 hours}
    \label{fig:load_solar_profile}
    \squeezeup
\end{figure}

The proposed control is tested on a large unbalanced three-phase 1500 node system based on the IEEE123 bus feeder \cite{noauthor_ieee_nodate-1}. To create a realistic simulation, the 123 bus system is further expanded with detailed secondary side house-load modeling at 120 volts resulting in 1500 nodes as shown in \figurename \ref{fig:test_ckt} using GridLAB-D platform, an open-source agent-based simulation framework for smart grids developed by Pacific Northwest National Lab \cite{chassin_gridlab-d:_2014}. Each residential load is modeled in detail with ZIP loads and temperature dependent HVAC load \cite{schneider_multi-state_2011}. Diversity and distribution of parameters within the residential loads is discussed in \cite{fuller_evaluation_2012}. The feeder is populated with 1280 residential houses with approximately 6 MW peak load.  Inverter ratings are considered 1.1 times the panel ratings. Uniformly distributed solar PVs throughout the feeder create lesser problems than the PV units distributed in one area of the feeder. Therefore, to demonstrate the effectiveness of the control in more severe case, PV units are distributed randomly at 500 houses only in right half of the feeder. Temperature and solar irradiance data for January 2, 2011 is obtained from publicly available NREL data for Hawaii \cite{sengupta_oahu_2010}. Load and solar profiles for the day have been shown in \figurename \ref{fig:load_solar_profile}.  {Voltage regulator at the substation is not de-activated and has a time delay of 5 minutes. It is not expected to interfere in proposed VVC due to difference in time-scale.}

\vspace{-0.2em}

\subsection{Performance Metrics}
We will be using three performance metrics to evaluate the proposed control approach. First metric is mean steady state error ($MSSE$) which evaluates control set-point tracking performance and
calculated as
\begin{equation}
\small
\label{eq:msse}
MSSE=\sum_{i=1}^n \sum_{t=1}^h \frac{|V_{t,i}-\mu|}{h}.\frac{1}{n}\times100
\end{equation}
where $n$ is the total number of solar PV units and $h$ is total time duration. A lower MSSE denotes better set point tracking performance of the control.
Second metric is flicker count ($FC$) where one flicker violation at one house is considered when $VF$ 
value, as defined in (\ref{eq:vf}), exceeds $VF_{lim}$. 
The total number of such flicker violations at all of the houses is termed as $FC$. Higher value of this metric is an indication of lesser power quality and an oscillatory voltage profile that in turn indicates the possibility of unstable control.
The third metric is voltage violation index ($VVI$) is the total number of voltage violations at all of the buses. Based on ANSI standards \cite{noauthor_ansi_2016}, 
a voltage violation is counted if the voltage at a bus violates either 1) 1.06-0.9 pu band instantaneously (range A) or 2) 1.05-0.95 pu band continuously for 5 minutes (range B).

\subsection{Results}
\begin{figure}
	\centering
	\includegraphics[trim=0in 0in 0in 0in,width=3.4in]{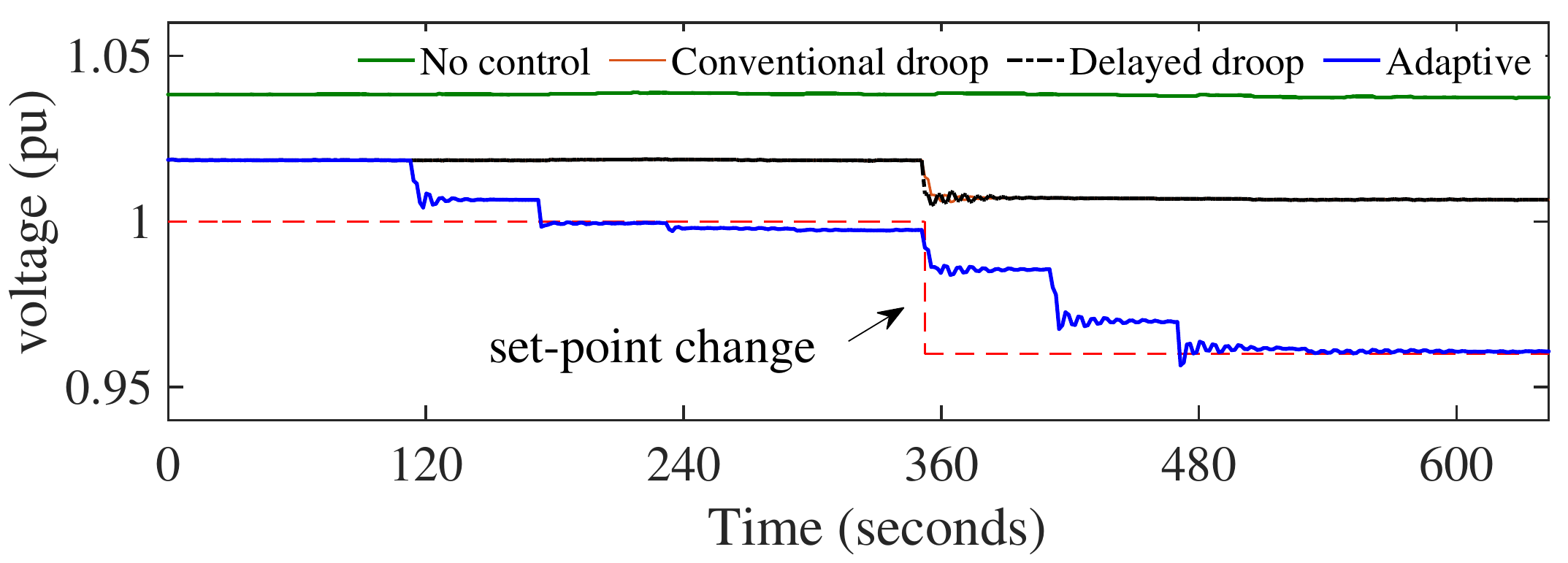}
    \vspace{-0.5em}
    \caption {Voltage profiles to compare set-point tracking performance of adaptive control with other methods}
    \label{fig:setpoint_change}
\end{figure}

In this section, we will demonstrate the effectiveness of the proposed control scheme in a wide range of external disturbances and operating conditions. The adaptive control performance (blue) will be compared with existing droop controllers i.e. conventional (orange) and delayed droop (black) as shown in \figurename\ref{fig:setpoint_change}. Dashed and solid red lines denote the voltage violation limits and voltage set-point respectively. Voltage profiles and parameters dispatched are shown at a randomly chosen solar PV unit at bus 92 whereas the performance metrics are calculated for the whole system as shown in Table \ref{tab:24_hours}. Outer loop horizon $T\!=\!1$ minute and $k^d\!=\!4$ are considered.
\subsubsection{Control Performance on Static Load Conditions}

In order to evaluate how well the proposed control tracks a set-point and minimizes SSE, a sudden set-point change is applied at the static load conditions in \figurename \ref{fig:setpoint_change}. Load and solar conditions at 11 AM are used for this purpose. 
Set point $\mu$ is changed from 1 to 0.96 pu at the middle of the simulation for all the inverters. As expected, the voltage profile in the no control case is steady. The conventional and delayed droop control have similar steady state performance which fails to track the set-point accurately and settles down with high SSE value. The adaptive control scheme, however, is able to track the set-point accurately by adapting $q_p$ parameter. It verifies the adaptive control's SSE minimization capability.

\subsubsection{Dynamic Tests with Daily Load and Solar Variation}
\begin{figure}
	\centering
	\includegraphics[trim=0in 0in 0in 0in,width=3.5in]{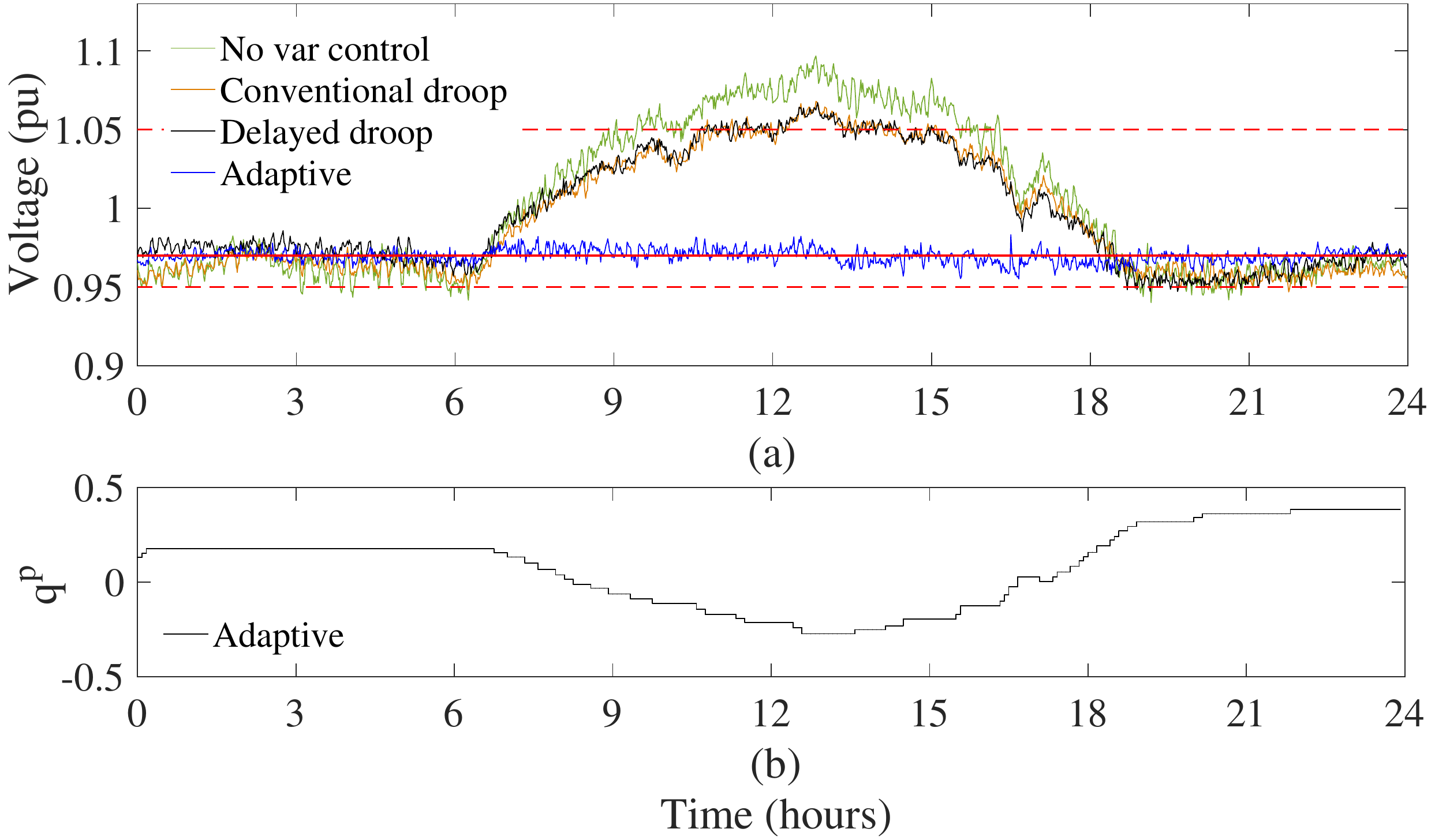}
    \vspace{-1.5em}
    \caption {Comparison of adaptive control performance throughout the day: a) voltage profile; b) dispatch of error adaptive parameter ($q^p$)}
    \label{fig:24hr_voltage}
    \vspace{-1.5mm}
\end{figure}

A day-long load and solar profile can be seen as continuous external disturbances in the system. \figurename\ref{fig:24hr_voltage}(a) shows that during the daytime, non-adaptive droop controls are not able to track the set point voltage which might lead to voltage violations e.g. around 12 noon when the solar generation is at peak. $\mu=0.97$ and homogeneous conservative settings $(m=3)$ are used for conventional and delayed control. Whereas the adaptive control  adapts it’s parameters at each bus differently to keep a flat voltage profile throughout the day; note, this may not be entirely desirable for the utility, or the owners, due to increased var flows, but rather indicates the flexibility of the system for applications such as CVR, loss minimization etc. {The OLTC tap counts are significantly reduced from 31 in no var control case to 5 in adaptive VVC case which improves the tap changer's life span. However, the tap changers can also be pre- dispatched based on the day-ahead load/solar profiles using a supervisory control, if required.} \figurename\ref{fig:24hr_voltage}(b) shows the dynamic dispatch of adaptive error parameter $q^p$ at bus 92. The performance metrics for the whole system are compared in Table \ref{tab:24_hours}. In this case, high MSSE in delayed control is because of selecting a conservative slope setting which can be improved by choosing higher slope, however, it will make the control highly vulnerable to sudden external disturbances as demonstrated in the next results. Whereas due to its decoupled functionality, the proposed control is capable of achieving near zero MSSE even at conservative settings, thus not making system prone to instability or voltage flicker.

\begin{table}
\caption{Performance Metrics Comparison for 24-hour profile}
\vspace{-1em}
\label{tab:24_hours}
\begin{center}
\renewcommand{\arraystretch}{1}
\begin{tabular}{|c|c|c|c|c|}
\hline
Metrics & \shortstack{No\\control} & \shortstack {Conventional \\droop} & \shortstack{Delayed  \\droop} & \shortstack{Adaptive \\control}\\
\hline
MSSE & 5.2\% & 4.3\% & 4.3\% & 0.3\%\\
VVI & $5\times 10^5$ & 21853 & 16137 & 0\\
\hline
\end{tabular}
\end{center}
\squeezeup
\end{table}

\subsubsection{Dynamic Tests with Sudden External Disturbances}
Reliable performance under external disturbances is a unique feature of the proposed control. To demonstrate it, the control is tested with sudden external disturbances. A smaller window of 1-2 hours is considered when solar is at its peak to observe the most severe impact of disturbances.

\begin{figure}
	\centering
	\includegraphics[trim=0in 0in 0in 0in,width=3.4in]{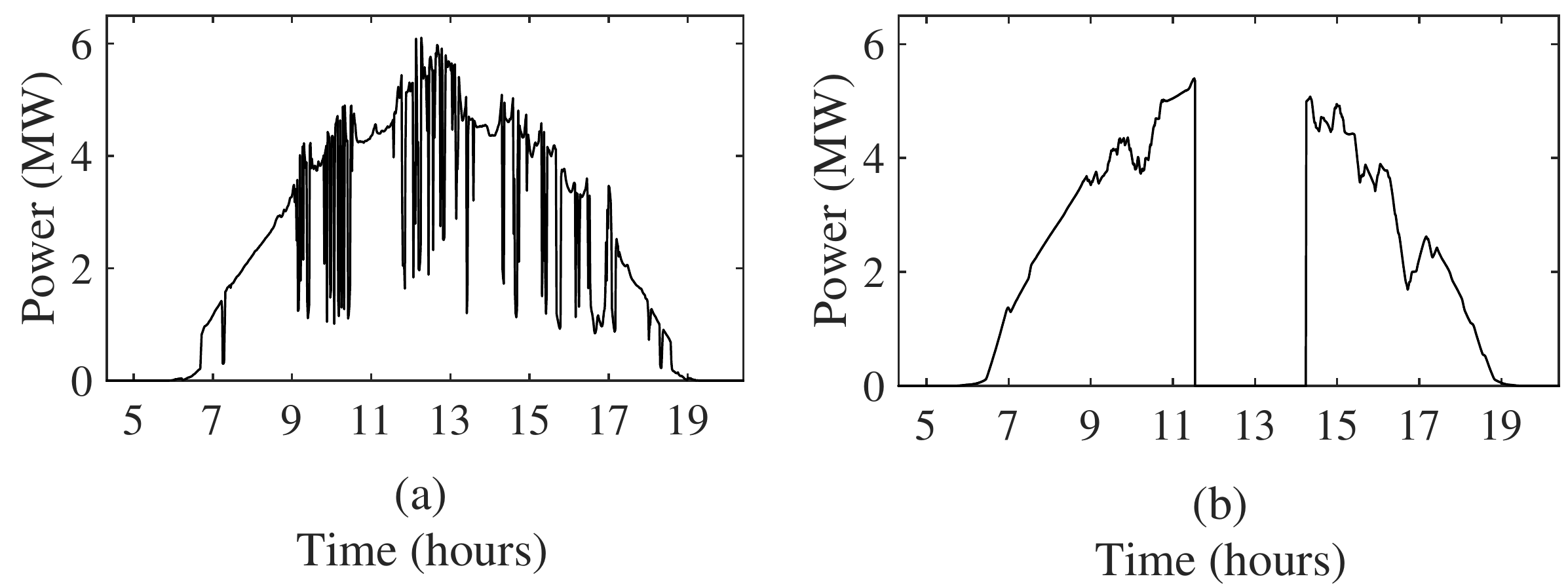}
    \vspace{-0.5em}
    \caption {Solar profile with a) cloud intermittency b) cloud cover}
    \label{fig:cloud_profile}
    \vspace{-0.1mm}
\end{figure}

\begin{figure}
	\centering
    	\includegraphics[trim=0.0in 0.0in 0in 0in,clip,width=3.6in]{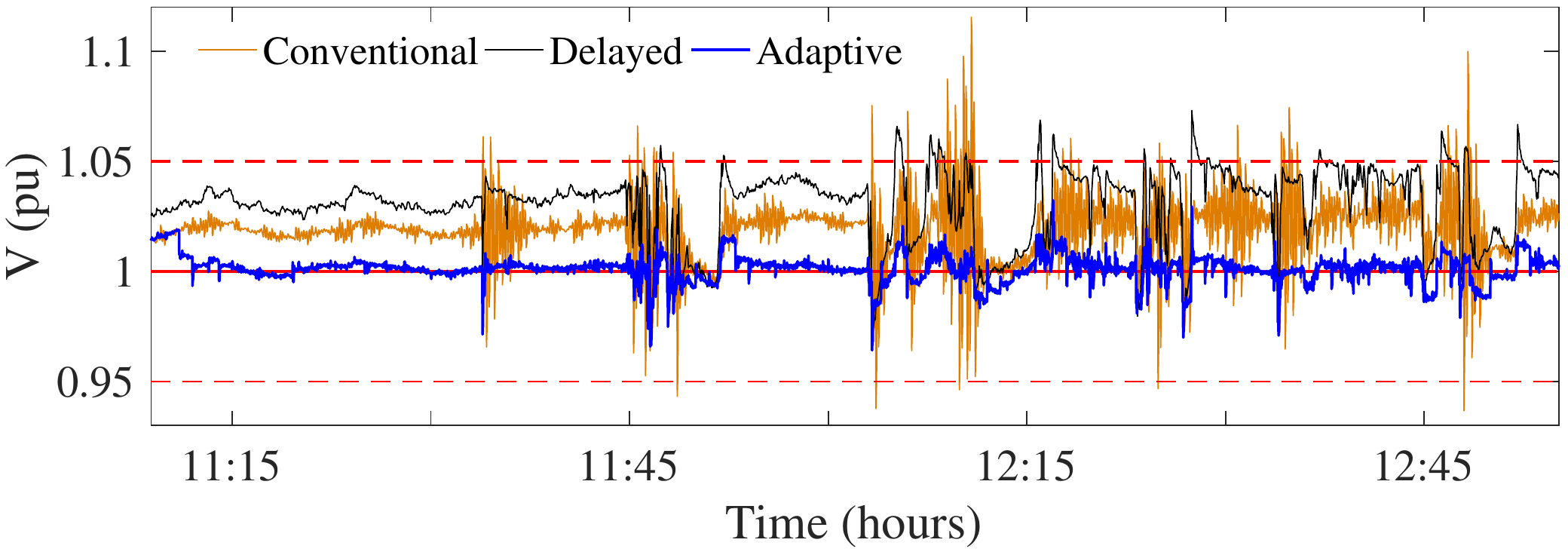}
    \vspace{-1.9em}
    \caption {Control performance comparison under cloud intermittency}
    \label{fig:cloud_intermittent_voltage}
     \vspace{0.5em}
     
    \includegraphics[trim=0in 0.3in 0in 0in,clip,width=3.55in]{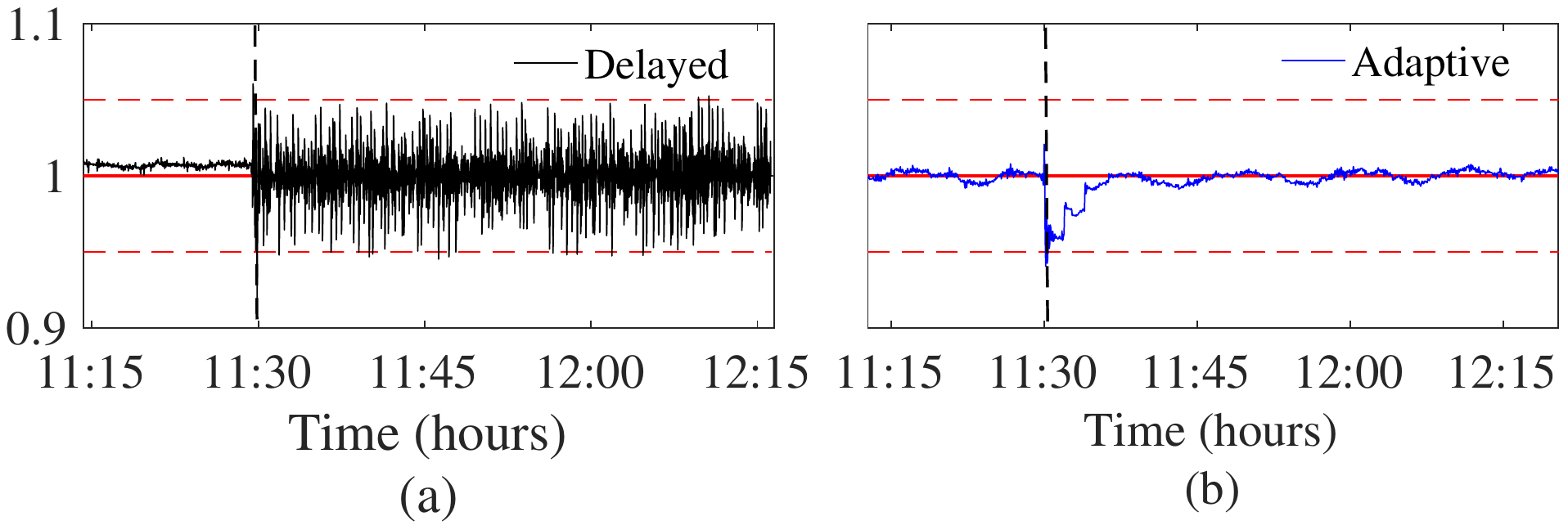}
    \vspace{-1.5em}
    \caption {Control performance comparison under sudden cloud cover}
    \label{fig:cloud_cover_voltage}
    \vspace{-0.8em}
\end{figure}

\textit{a) Sudden cloud cover and cloud intermittency}: Usually cloud covers cause two types of disturbances in PV generation i.e. intermittency and sudden drop in the generation as shown in \figurename \ref{fig:cloud_profile}(a) and (b) respectively. Cloud intermittency data of 30 seconds scale is considered. Set-point $\mu=1$ and $m=5$ are used for the non-adaptive controls. \figurename \ref{fig:cloud_intermittent_voltage} shows how cloud intermittency causes high voltage fluctuations in conventional control which leads to violations. Delayed control reduces the flicker significantly compared to conventional (from 6919 to 107), however, still results in a good number of violations due to high SSE as shown in Table \ref{tab:cloud_intermittent}. Though, the effect of intermittency is visible in adaptive control profile (\figurename \ref{fig:cloud_intermittent_voltage}), it manages to achieve zero indices of flicker and violations. It shows the effectiveness of control in faster disturbances.

On the other hand, using non-conservative settings $(m=10)$ to decrease violations can cause stability issues with sudden cloud cover as shown in \figurename \ref{fig:cloud_cover_voltage}. At 11.30 AM, a cloud cover results in a sudden drop in real power generation (\figurename \ref{fig:cloud_profile}(b)) which frees the inverter capacity. Since conventional and delayed controls utilize all the free capacity immediately without monitoring, it increases the slope by a significant amount and results in voltage oscillations as shown in  \figurename \ref{fig:cloud_cover_voltage}. Whereas, the adaptive control dynamically regulates the  settings in real-time to ensure stable voltage profile as well as quick restore of the set-point tracking. fhl{Moreover, the proposed control can also be integrated with other smart grid applications where sudden change in real power generation is experienced such as when PV inverters are providing virtual inertia to the system.} 

\begin{table}
\vspace{-0.5em}
\caption{Performance Metrics Comparison For Intermittent Solar-Profile For A Two-Hour Window}
\vspace{-1em}
\label{tab:cloud_intermittent}
\begin{center}
\renewcommand{\arraystretch}{1.2}
\begin{tabular}{|c|c|c|c|c|}
\hline
Metrics & No control & Conventional & Delayed  & Adaptive\\
\hline
MSSE & 3.5\% & 2.00\% & 2.00\% & 0.40\%\\
VVI & $5\times 10^5$ & 21853 & 16137 & 0\\
FC & 122 & 6919 & 107 & 0\\
\hline
\end{tabular}
\end{center}
\end{table}

\begin{figure}
	\centering
	\includegraphics[trim=0in 0.3in 0in 0in,clip,width=3.55in]{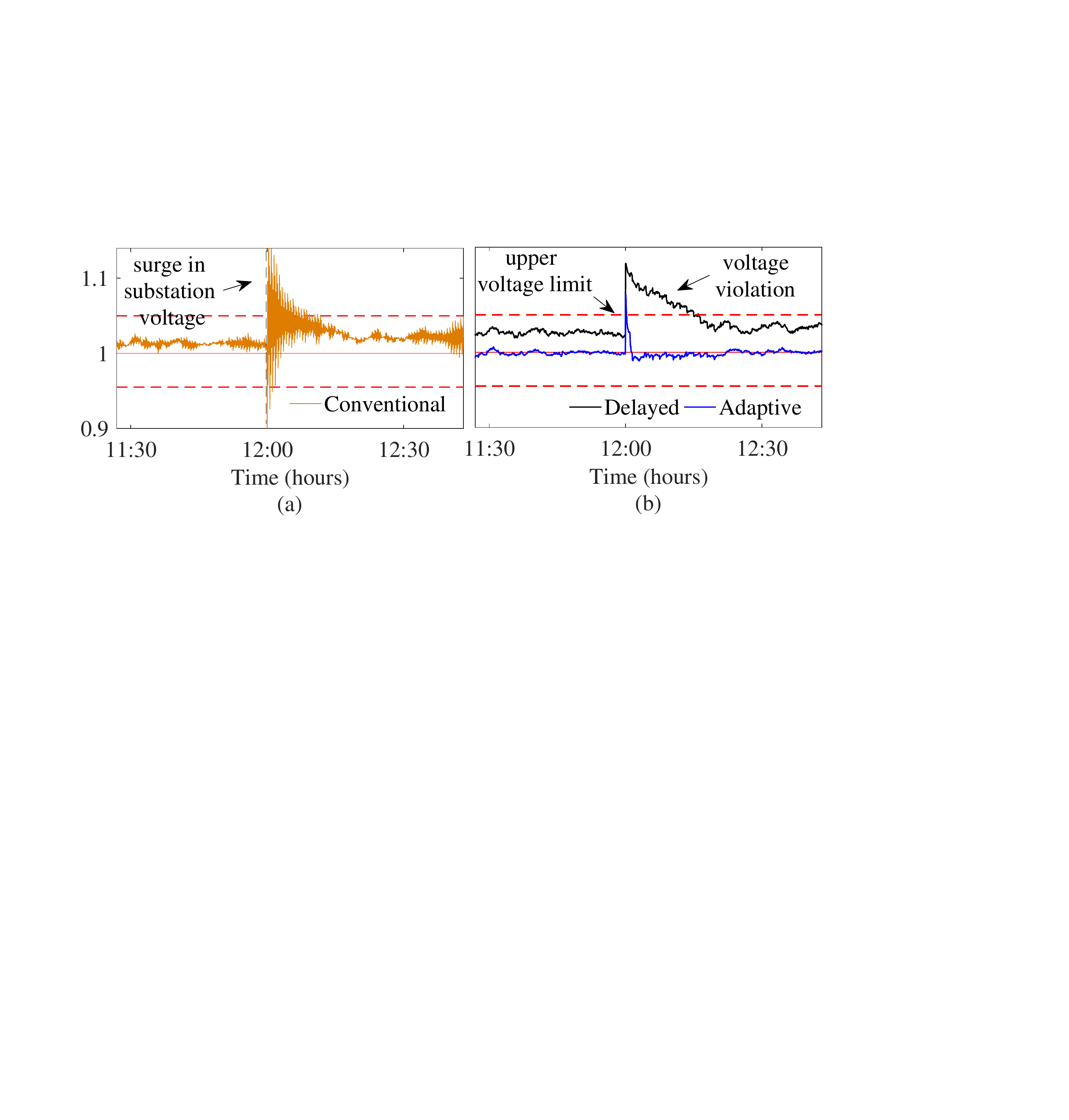}
    \vspace{-1.5em}
    \caption {Impact of change in substation voltage on control performances}
    \label{fig:substation_change}
    \vspace{-.5em}
\end{figure}

\textit{b) Change in substation voltage}: The primary side of substation voltage keeps changing due to changes in the transmission systems. Conservative setting $(m=5)$ is used here for non-adaptive controls. In  \figurename \ref{fig:substation_change}, at 12 noon, the feeder experiences a surge in primary substation voltage from 1 to 1.07 pu. Conventional control experiences high voltage oscillations. Delayed control does not experience voltage flicker but since it cannot reduce the SSE on its own, it waits for substation tap changer to operate to bring voltage within the limit again. Whereas adaptive control suffers from few instantaneous violations but immediately starts re-tracking the set-point, thus avoiding violations for long time period.

Similarly, the proposed approach can also be verified in several other scenarios such as sudden load drop, increasing PV penetration etc. {However, in order to deploy the PV inverter VVC effectively, the customers need to be given proper incentives by the utilities to allow thier inverters to participate in var support. Moreover, in emergency situations, the real power curtailment might be necessary for which the owners needs to be compensated. Therefore, some utilities are moving towards the utility owned solar and encouraging community solar projects which are more viable in terms of providing VVC benefits.}

\vspace{-1em}
\section{Conclusion}
\vspace{-1mm}
In this study, a real-time adaptive and local VVC scheme with high PV penetration is proposed to addresses two major issues of conventional droop VVC methods. First,  the proposed adaptive droop framework enables VVC to achieve high set-point tracking accuracy and control stability (low voltage flicker) simultaneously without compromising either. Second, it enables dynamic self-adaption of control parameters in real-time 
under wide range of operating conditions/external disturbances. All this is achieved while keeping the control 
compatible with the integration standards (IEEE1547) and utility practices (Rule 21). 
The satisfactory performance is demonstrated by comparing with existing droop methods in several cases on a large unbalanced distribution system.

Being local in nature, the proposed control might be less effective for system-wide optimization, however, the proposed VVC framework can easily be combined with centralized approaches. In fact, due to its tight voltage regulation feature and adaptive nature under external disturbances, it facilitates the use of PV inverters for other system-wide volt/var applications such as CVR, loss minimization, var support to the grid, etc. The integration with supervisory control and coordination with conventional regulators will be explored in the future work.

\ifCLASSOPTIONcaptionsoff
  \newpage
\fi


%


\bibliographystyle{IEEEtran}
\bibliography{IEEEabrv,References_short}

%
 \begin{IEEEbiography}
 [{\includegraphics[width=1.1in,height=1.2in,clip,keepaspectratio]{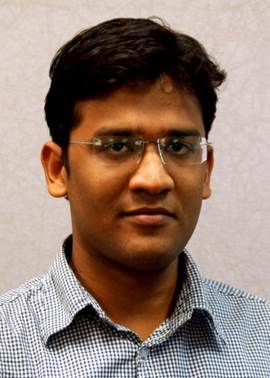}}]{Ankit Singhal}
  (S'13) received the B.Tech. degree in electrical engineering
from the Indian Institute of
Technology-Delhi, India. He is currently
a Ph.D. student in the Department of Electrical and Computer Engineering at Iowa State
University, Ames, IA, USA.

His research interests include renewable integration, impact of high PV penetration on distribution and transmission network, smart inverter volt/var control and transmission-distribution co-simulation.
 \end{IEEEbiography}
 \begin{IEEEbiography}
 [{\includegraphics[width=1in,height=1.2in,clip,keepaspectratio]{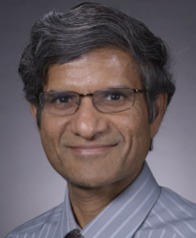}}]{Venkataramana Ajjarapu}
  (S'86-M'86-SM'91-
F'07) received the Ph.D. degree in electrical engineering from the University of Waterloo, Waterloo, ON, Canada, in 1986.

Currently, he is a Professor in the Department of Electrical and Computer Engineering at Iowa State University, Ames, IA, USA. His present research is in the area of integration of distributed energy resources, reactive power planning, voltage stability analysis, and nonlinear voltage phenomena.
 \end{IEEEbiography}
 
  \begin{IEEEbiography}
 [{\includegraphics[width=1in,height=1.2in,clip,keepaspectratio]{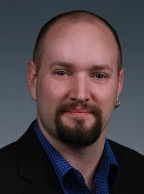}}]{Jason Fuller}
  (S'08-M'10) received the B.S. degree
in physics from the University of Washington,
Seattle, and the M.S. degree in electric engineering from Washington State University, Pullman.

He is currently the manager of Electricity Infrastructure group at the Pacific Northwest National Laboratory. His main areas of interest are distribution system analysis and renewable integration. He is currently the Secretary of the Distribution System Analysis Subcommittee’s Test Feeder Working Group.
\end{IEEEbiography}

  \begin{IEEEbiography}
 [{\includegraphics[width=1in,height=1.2in,clip,keepaspectratio]{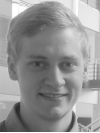}}]{Jacob Hansen}
 (S'12-M'14) received
the B.Sc. and M.Sc. degrees in electronics
engineering and information technology from
Aalborg University, Aalborg, Denmark.

Currently, he is employed at Pacific Northwest National Laboratory as a Power System Control Engineer. He was a Visiting Student with the Active Adaptive Control Laboratory, Massachusetts Institute of Technology, Cambridge, MA, USA,
from 2013 to 2014. His current research interests include smart grid, power systems, decentralized
control, and control systems in general.
\end{IEEEbiography}







\end{document}